\documentclass{amsart}
\usepackage{amsmath,amsthm,amscd}

\newtheorem{thm}{Theorem}[subsection]
\newtheorem{conj}{Conjecture}[subsection]
\newtheorem{cor}[thm]{Corollary}
\newtheorem{lem}[thm]{Lemma}
\newtheorem{prop}[thm]{Proposition}

\theoremstyle{definition}
\newtheorem{defn}[thm]{Definition}

\newtheorem{spec-case}[thm]{Special Case}

\theoremstyle{remark}
\newtheorem{rem}[thm]{Remark}

\newtheorem{example}[thm]{Example}
\newtheorem{claim}[thm]{Claim}

\newtheorem{question}[thm]{Question}

\newcommand{\cF}{\mathcal{F}}

\newcommand{\cK}{\mathcal{K}}

\newcommand{\cM}{\mathcal{M}}
\newcommand{\cN}{\mathcal{N}}
\newcommand{\cO}{\mathcal{O}}
\newcommand{\cR}{\mathcal{R}}

\newcommand{\fa}{\mathfrak{a}}

\newcommand{\fm}{\mathfrak{m}}

\newcommand{\fp}{\mathfrak{p}}

\newcommand{\fM}{\mathfrak{M}}
\newcommand{\fN}{\mathfrak{N}}

\newcommand{\bbI}{\mathbb{I}}
\newcommand{\bbN}{\mathbb{N}}
\newcommand{\NN}{\mathbb{N}}
\newcommand{\bN}{\mathbb{N}}
\newcommand{\bbP}{\mathbb{P}}

\newcommand{\bbZ}{\mathbb{Z}}

\newcommand{\bbQ}{\mathbb{Q}}
\newcommand{\QQ}{\mathbb{Q}}
\newcommand{\bbR}{\mathbb{R}}
\newcommand{\bbC}{\mathbb{C}}
\newcommand{\CC}{\mathbb{C}}
\newcommand{\bbV}{\mathbb{V}}

\newcommand{\bR}{\mathbf{R}}

\newcommand{\gl}{\lambda}

\newcommand{\lrarrow}{\longrightarrow}
\newcommand{\la}{\longrightarrow}

\DeclareMathOperator{\SHom}{\mathcal H {\mathit o \mathit m}}
\DeclareMathOperator{\SExt}{\mathcal E {\mathit x \mathit t}}

\DeclareMathOperator{\Hom}{Hom}\DeclareMathOperator{\Ext}{Ext}

\DeclareMathOperator{\ann}{Ann}
\DeclareMathOperator{\Ann}{Ann}

\DeclareMathOperator{\core}{core}
\DeclareMathOperator{\grcore}{graded core}

\DeclareMathOperator{\adj}{adj}
\DeclareMathOperator{\codim}{codim}

\DeclareMathOperator{\Min}{Min}
\DeclareMathOperator{\Proj}{Proj}
\DeclareMathOperator{\Spec}{Spec}
\DeclareMathOperator{\spec}{Spec}
\DeclareMathOperator{\proj}{Proj}
\DeclareMathOperator{\gr}{gr}

\title[]{On a Non-Vanishing Conjecture of Kawamata and
 the Core of an Ideal } 
\author{Eero Hyry}
\address{Department of Mathematics\\
PL 4 (Yliopistonkatu 5), FIN-00014 HELSINGIN YLIOPISTO, Finland}
\email{Eero.Hyry@helsinki.fi}
\thanks{Research of Hyry supported by the National Academy of Finland,
project number 48556.}
\author{Karen E. Smith}
\address{Department of Mathematics\\ University of Michigan, Ann Arbor, MI, 48109-1109}
\email{kesmith@umich.edu}
\thanks{Research of Smith partially  
 supported by a US Fulbright fellowship and 
 the National Science Foundation DMS 00-70722. Both authors are
grateful for the hospitality of Jyv\"akyl\"a University where the
 bulk of this work was carried out.}
\subjclass{Primary  13 ; Secondary 14 .}
\keywords{effective non-vanishing, adjoint, core, multiplier ideal,
Shokurov's non-vanishing}
\begin{document}
\begin{abstract} 
We show, under suitable hypothesis which are sharp in a certain sense,
that the core of an $\fm$-primary  ideal in a regular local ring of dimension 
$d$ is equal
to the adjoint (or multiplier) ideal of its $d$-th power.
This generalizes the fundamental
formula for the core 
of an integrally closed ideal in a  two dimensional regular local ring 
due to Huneke and Swanson.
We also find a generalization of this 
result to singular (non-regular)
settings, which we show to be intimately related to the 
problem of finding non-zero sections of ample line bundles on projective
varieties. In particular, we show that a graded analog of our 
formula for core would imply a 
remarkable conjecture of Kawamata predicting that every adjoint ample
line bundle on a smooth variety admits a non-zero section. 
\end{abstract}

\maketitle

\tableofcontents

\section{Introduction}\label{sec}

Let $I$ be an  ideal in a commutative ring. By definition,
the {\it core} of $I$ is the intersection 
of all sub-ideals having the same  integral closure as $I$.
Because the notion of integral closure 
is so fundamental, the core is a natural and interesting object. Originally 
 defined by Rees and Sally in \cite{RS}, the first substantial 
progress understanding core is due to Huneke and Swanson in \cite{HS},
who proved that   the core of an integrally closed $\fm$-primary
ideal $I$ in a two dimensional regular local ring 
 is equal to the adjoint (or multiplier) ideal of $I^2$.
Since then, the algebraic properties of core have been 
thoroughly studied;
  see, for example, \cite{CPU} and \cite{CPU2}.
Our own interest in the core is  motivated
by
seemingly unrelated  geometric concerns:  
the core of a certain ideal governs whether or not an ample
line bundle on a projective algebraic variety has a non-zero global 
section.

The purpose of this paper is two-fold. On the algebraic side, 
we find structure theorems  for the core 
naturally generalizing  the  Huneke-Swanson results to  the 
 higher dimensional
and singular case.
On the geometric side, we show how a sufficiently good understanding of
the core of a certain ideal in a very special kind of graded ring
would settle a remarkable
 conjecture
of Kawamata
 predicting that ``adjoint'' nef divisors are always effective.
In particular, 
our higher dimensional 
singular version of the Huneke-Swanson formula for core 
 can be viewed as 
 a local analog of Kawamata's conjecture.
Our work shows that  
 commutative algebraists and algebraic geometers,
working completely independently of each other and motivated by
very different problems,  both discovered and conjectured 
different facets of the same beautiful---and still largely buried--- 
mathematical 
diamond. 

It is perhaps no surprise that multiplier (adjoint) ideals arise in
the search for 
 non-vanishing theorems
 for nef line bundles.
In recent years, multiplier ideals
 have found many rich applications to algebraic
geometry, particularly  
 to issues of effectiveness; see \cite{Dem}, \cite{Siu}. 
Originally defined as ideals of holomorphic functions belonging to a
certain weighted $L^2$-space, 
 multiplier ideals
can also be developed  purely algebraically (see \cite{L3}) 
or algebro-geometrically (see \cite{Ein} or \cite{laz}). 
Recently, multiplier ideals have been used to prove surprising
new results in 
commutative 
algebra as well; see 
\cite{ELS1} and 
 \cite{ELS2}.

\bigskip
\subsection {The motivating Geometry}
Let $X$ be a complex smooth projective variety and let $L$ be an ample
line bundle on $X$. Being ample,
 the line bundle 
$$
 L^n = \underbrace{L\otimes L \otimes \dots \otimes L}_{\text{$n$ times}}
$$
has many global sections for large enough $n$.
On the other hand,  $L$ itself need not have any 
sections at all. 
A fundamental unanswered question in
algebraic geometry is this: what $n$ is large enough so that $L^n$ has even one
non-trivial 
global section?

As stated, there is no general answer to this question; 
there is no uniform $n$ that works for all line bundles on a given variety
$X$. So we must restrict attention to particular classes of bundles. 
For example, 
suppose that  $L$ is an adjoint bundle, that is, that 
$L$ is of the form 
$K_X \otimes H$  where $H$ is some ample line bundle and 
   $K_X$ is
the canonical bundle   of $X$.
The Kodaira vanishing theorem implies that the higher cohomology groups 
of $L$ must all vanish, but as pointed out by Ambro in \cite{Am},
not a single
example is known  in which  the  zeroth cohomology group also vanishes. 
In fact, even for  numerically effective (``nef'') line bundles $L$ of the
form 
$K_X \otimes H$ where $H$ is ample, the celebrated theorem of 
 Shokurov
guarantees that $L^n$ has many sections for 
sufficiently large $n$  \cite{Shok};  again, in all known examples,
$L$ itself already has a non-zero section.

\smallskip
In \cite{kaw1}, 
Kawamata  conjectured that  every numerically 
effective line bundle $L$ adjoint  to an  ample line
bundle
  must
have a non-zero global section: that is, if $L$ is of the form 
$K_X \otimes H$ for some ample $H$, then
 $L^n$ has a non-zero section for {\it every}
$n \geq 1$.
More generally, 
 Kawamata stated his  effective non-vanishing 
conjecture in a singular, logarithmic  
 version, indeed, under the same  general 
hypothesis of Shokurov's  
non-vanishing theorem (or of the ``Base Point Free Theorem'') 
\cite{kaw1}:
\begin{conj}\label{mainkaw}
Let $D$ be any numerically effective (``nef'')
 Cartier divisor on a normal projective variety $X$.
 If there exists an effective $\bbR$-divisor $B$ 
such that the pair $(X, B)$ is  Kawamata log terminal
{\footnote{One can get  interesting, 
 but less technical statements by taking
 $B$ to be zero and  $X$ to be Gorenstein:
then the condition that $(X, B)$ is Kawamata log terminal amounts to
$X$ having rational singularities.
For the basic terminology of singularities of pairs used here,
we recommend \cite{Ko} or \cite{KM} as a general reference.}}
 and such that the $\bbR$-Cartier divisor 
$D - (K_X + B)$ is big and nef, then the line bundle $\cO_X(D)$
has a non-zero global section. 
\end{conj}

This remarkably strong conjecture was first raised as a 
question by Ambro in \cite{Am}.
For curves it is trivially true and for surfaces it can be shown
using a Riemann-Roch argument (\cite[Theorem 3.8]{kaw1}), but in higher 
dimensions it seems quite surprising.
For example,
Conjecture \ref{mainkaw}  predicts that
  the linear system of {\it every} nef divisor on a smooth Fano variety 
is non-empty; in particular, every Fano variety admits an effective 
anti-canonical divisor.
For smooth varieties for which $K_X$ is trivial,
Conjecture  \ref{mainkaw} asserts that every big and nef line bundle 
has a non-zero section; 
in particular, 
 every ample line bundle on a Calabi-Yau
manifold admits a non-trivial global section. 
At the opposite extreme,  
 Conjecture \ref{mainkaw}  guarantees the existence of non-zero sections
for  the bi-canonical bundle on a smooth 
 minimal model of general type.

\medskip

\subsection{From Geometry to Algebra.}
 Kawamata has shown that  in order to prove Conjecture \ref{mainkaw},
it suffices to 
 prove the case where $D$ is ample \cite[Theorem 2.2]{kaw1}.
This opens up the possibility of using commutative algebra. 

 Given a pair $(X, D)$ consisting of
a normal projective variety $X$ and an ample divisor  $D$,
  consider the section ring
$$
S = \bigoplus_{n \in \bbN} H^0(X, \cO_X(nD)),
$$
a normal, finitely generated  graded ring 
whose associated projective scheme recovers $X$. 
As we describe in detail in Section  6,
Kawamata's Conjecture leads naturally to an equivalent
 statement about the 
 {\it graded  core} of a certain submodule of
the canonical module  $\omega_S$ of $S$.
This statement takes its simplest form in the case where 
 $D = -K_X$ is the anti-canonical bundle on a smooth Fano variety.
In this case, the section ring $S$ is Gorenstein and 
the line bundle $-K_X$ will have a non-zero 
section if the formula
\begin{equation}\label{conj}
\grcore{(I)} = \adj{(I^{d})}, \, \, {\text { where }} d = \dim S,
\end{equation}
holds 
for the integrally closed ideal $I = S_{\geq n}$ generated by elements of
degrees at least $n$, for some (equivalently, every) fixed large $n$.
Here the graded core of $I$ is the intersection of
all homogeneous subideals whose integral closure is   $I$,
and $ \adj{(I^{d})}$ is the adjoint ideal (or multiplier ideal) 
of $I^{d}$. (The definition is recalled in Section \ref{adjoint}).

For more general pairs $(X, D)$,  the section ring $S$ will not be Gorenstein,
so the usual notion of a multiplier ideal need not be defined. 
However, we overcome this difficulty by
 using an alternate construction of {\it adjoint module};
this adjoint module is  a submodule of $\omega_S$ rather than an ideal of $S$,
agrees with the adjoint ideal for local Gorenstein rings after 
fixing an isomorphism $\omega_S \cong S$, 
but seems better suited for working on singular varieties.{\footnote{The
point is that its definition does not require the existence of a relative
canonical divisor.}}
In this more general setting, we show that 
  $D$ has a section if 
\begin{equation}\label{conj2}
\grcore{(I\omega_S) } = \adj{(I^{d}\omega_S)},
 \, \, {\text { where }} d = \dim S,
\end{equation}
holds for the ideal $I = S_{\geq n}$, $n \gg 0$,  where 
$\grcore{(I\omega_S) }$ is the intersection of the submodules $J\omega_S$ 
in $\omega_S$, as $J$ ranges  over all {\it  homogeneous} ideals 
$J$ whose integral closure is $I$,
and
$\adj{(I^{d}\omega_S)}$ is the adjoint module, as defined in 
  Remark 
\ref{mult.ideal}.

Thus it may be possible to prove 
 Kawamata's Conjecture 
by 
 proving formulas  
(\ref{conj}) or (\ref{conj2}) for section rings
of divisors $D$ satisfying the conditions of 
Kawamata's hypothesis. Such section rings
are very special:
for example, in the Fano case above, the Rees ring $S[It]$
turns out to be Cohen-Macaulay and even to have rational singularities;
 even for more general pairs, the Rees ring is $S[It]$ still has the
very strong property that its irrelevant ideal is a Cohen-Macaulay module.
Under the strong conditions imposed by the geometric hypothesis, we 
hope to eventually prove formula (\ref{conj2})
and hence Kawamata's Conjecture.
What we do in the current paper is a local 
version of just that.

\medskip
\subsection{The Main Algebraic Results.}
Guided by the conjectural formulas (\ref{conj}) and (\ref{conj2}),
 we prove formulas for the core of an 
ideal $I$ in a local ring $A$ satisfying
 hypotheses  satisfied by the ideal $I = S_{\geq n} $ in the special
 section rings 
$S$ arising from a pair $(X, D)$ satisfying the assumptions
of Conjecture \ref{mainkaw}.
 For example, a special case of our main theorem produces
the following higher dimensional (and singular) version of the
Huneke-Swanson formula for core mentioned in the opening paragraph:
\begin{thm}\label{hun-swan}
Let $(A, \fm)$ be a Gorenstein local ring of dimension $d$
essentially of finite type over a field of characteristic zero,
 and let  $I$ be an $\fm$-primary ideal
of $A$. If the Rees ring $A[It]$ of $I$ 
has rational singularities,
then 
\begin{equation}\label{s}
\core (I) = \adj (I^{d}),
\end{equation}
where  $\adj (I^{d})$ denotes the
adjoint (or multiplier) ideal of the ideal $I^{d}$.
\end{thm}

Although the assumption above that the Rees ring has rational singularities
is quite strong, it is quite natural for two reasons.
First, 
it is satisfied when $S$ is the section ring of the anti-canonical
divisor on a Fano variety and $I = S_{\geq n}$, $n \gg 0$.
Thus 
 Theorem \ref{hun-swan} can
be viewed as a ``local version'' of  (\ref{conj})
and hence a local version 
of  (a special case of)  Kawamata's Conjecture.
Second, it
always holds
whenever $I$ is an integrally closed ideal in a two-dimensional
regular local ring, the setting of the Huneke-Swanson formula.
In fact, we also show that under some further restrictions,
rational singularities of $A[It]$ is necessary and sufficient for
 formula (\ref{s}) to hold
 (see Theorem \ref{ratsingcor}).

\smallskip
Passing away from Gorenstein (and even Cohen-Macaulay) case,
as we must in order to prove Kawamata's
Conjecture in full generality,
 we arrive at a similar result inside the canonical module:

\begin{thm}
\label{lockaw}
Let $(A, \fm)$ be a local ring essentially of finite type over 
a field of characteristic zero, and let $I$ be an $\fm$-primary ideal
such that the irrelevant ideal of the Rees ring $A[It]$ is a Cohen-Macaulay
$A[It]$-module.
Then if $Y = \proj A[It]$ has rational singularities, 
$$
\core{(I\omega_A)} = \adj (I^d\omega_A) := \Gamma(Y, I^d\omega_{Y})
$$
as submodules of $\omega_A$ where $d$ is the dimension of $A$.
\end{thm}

Theorem \ref{lockaw}, then, can be interpreted as a local version of
(\ref{conj2}), and so a local analog 
Kawamata's Conjecture.  Indeed, its hypothesis 
is precisely satisfied by 
  a section ring $S$ of a pair $(X, D)$ satisfying the hypothesis
of Conjecture \ref{mainkaw}. Therefore, for such a section ring, Theorem
\ref{lockaw} implies that 
\begin{equation}\label{s1}
\core(I\omega_S) = \adj(I^d\omega_S)
\end{equation}
where $ I = S_{\geq n}$ is the ideal generated by the
elements of degrees at least $n \gg 0 $ and $d$ is the dimension of $S$. 
We conjecture that furthermore
\begin{equation}\label{f}
\grcore(I\omega_S) = \adj(I^d\omega_S)
\end{equation}
holds, which, as  we show in Section 6, 
 implies 
Conjecture \ref{mainkaw}.
The difference between (\ref{s1}), which we prove to be true,
 and  (\ref{f}), which implies Kawamata's Conjecture, is that in the latter
we are intersecting over only homogeneous ideals $J$. On the other hand,
since both intersections are actually finite 
and there are  plenty of homogeneous reductions, 
it is reasonable to expect that perhaps equality holds. This, however,
seems to be a very subtle question.

\bigskip

Our formulas for core also elucidate the relationship between the core
and the  coefficient ideals of Aberbach and Huneke (see Corollary
 \ref{coeffcorecor}),
and  support some very general conjectures of
 Corso, Polini and 
Ulrich, which built on the Huneke-Swanson work; see \cite{CPU2}.
For example, we prove:
\begin{thm}\label{maincor}
Let $(A, \fm)$ be a Cohen-Macaulay local ring containing the
rational numbers and let $I$ be any equimultiple ideal of positive height.
If the Rees algebra $A[It]$ is Cohen-Macaulay, 
then
$$\core{I} = J^{r+1}: I^r \,\,\, {\text{ for }} r \gg 0,
$$
where $J$ is any minimal reduction of $I$.
\end{thm}
This formula for core is conjectured in more general settings
in \cite{CPU2}.

\bigskip

The format of the paper is as follows. 
Section 2 summarizes  some  background material and conventions,
while  recording various technical lemmas for later use.
Section 3 forms the technical heart of this paper.
Here we prove our main technical theorem, Theorem \ref{mainalg},
 describing the
core module $\core  (I\omega_A)$ of an  ideal in the canonical module
as an ``adjoint-type'' module.  
Theorem \ref{mainalg} requires a certain Brian\c con-Skoda-type
hypothesis,   and in Section 4,
we identify some natural classes of rings and ideals for which these
hypotheses are satisfied.
In Section 5, we  pull together these
results to deduce our main local  results,
including Theorems \ref{hun-swan} and \ref{maincor}.
Finally, in Section 6, 
we show how 
 Kawamata's Conjecture reduces to
a purely algebraic statement about the graded core analogous to our main 
results in the local case.

Throughout this paper, 
all rings and schemes are assumed Noetherian, and are
assumed to possess a dualizing complex.
The notation $(A, \fm)$ denotes a local ring whose unique maximal ideal is
named $\fm$. 
For an affine
scheme $\spec A$, we frequently abuse terminology by deliberating blurring the
difference between a quasi-coherent sheaf of modules on $\spec A$
 and the corresponding  $A$-module of its global sections.

\section{Algebraic Preliminaries.}\label{prelim}

This section summarizes some definitions,
 tools and conventions we will use.
 With the exception of the definition of core and graded core
in Subsection 1 and a few  technical results 
(likely to be well-known to experts),
 nearly all of this material can be unearthed from the
sources \cite{HIO}, \cite{RD}, and \cite{L2}. 
Readers
may prefer to 
skip this section and refer back only as necessary.

\subsection{Integral Closure, Reductions, and the Core}

Let $I$ be an ideal in a Noetherian ring $A$.
The {\it integral closure} of $I$ is defined as the set of
all elements $z$
in $ A$ satisfying a polynomial equation
$$
z^n + a_1 z^{n-1} + \dots + a_n = 0, 
$$
where $a_i \in I^i$. The integral closure $\overline I$
 of $I$ is an ideal of $A$ containing $I$ and 
contained in the radical of $I$.  
In the case where  $I$ is a homogeneous ideal of a graded ring, 
$\overline I$ is also homogeneous.

A {\it reduction}
 of $I$ is any sub-ideal of  $I$ having the same integral closure.
Equivalently, $J$ is a reduction of $I$ if there is a positive integer 
$r$ such that 
$$
I^{r+1} = JI^r.
$$
The smallest such $r$ is called the {\it  reduction number } 
of the pair $(I, J)$.
Equivalently, a sub-ideal $J$ is a reduction if and only if the
corresponding inclusion of Rees Rings
$$ 
A\oplus J \oplus J^2 \dots   \subset A\oplus I \oplus I^2 \dots
$$ 
is integral.

\smallskip
A {\it minimal reduction} of an ideal $I$ is a reduction that does
not properly contain any smaller reduction of $I$. Every ideal admits
minimal reductions. 
 If $I$ is an 
$\fm$-primary ideal in a local ring $(A, \fm)$ of dimension $d$, then 
any reduction generated by $d$ elements is a minimal reduction of $I$.
Conversely, if the residue field $A/\fm$ is infinite, then every minimal 
reduction is generated by $d$ elements. 

\smallskip
See \cite[\S4, \S10]{HIO} for more information on 
integral closures and reductions. 

\smallskip
\begin{defn}
The  {\it core} of an ideal is the intersection of all its reductions.
If the ideal is a homogeneous ideal in a graded ring, then its 
{\it graded core }
 is the intersection of all its homogeneous reductions. 
\end{defn}

\begin{rem}
Obviously, $\core{I} \subset \grcore(I)$ for any homogeneous
ideal in  a graded ring, but the inclusion can be strict in general,
even when $I$ has many homogeneous reductions; see Remark \ref{gradnotcore}.
\end{rem}

\smallskip
Our main interest from the point of view of Kawamata's 
conjecture 
is the ideal of elements of degrees at least $n$ in
a graded ring.  In this case we have the following simple descriptions of
integral closure, reductions, and the graded core. 

\begin{prop}\label{homred}
Let $S$ be an $\bbN$-graded reduced ring finitely generated 
 over an infinite field $S_0 = k$.
Suppose that the set of elements of $S$ of degree $n$ generate a cofinite
ideal $I$ (that is, that
 $S/I$ is finite dimensional over $k$).
Then
\begin{enumerate}
\item  $\overline I = S_{\geq n}$, the ideal generated by all elements of
degrees at least $n$. 
 \item  Every minimal homogeneous 
reduction of $I$ is generated by a 
system of parameters for $S$ consisting of elements of degree $n$. 
\item The graded core of $\overline{I}$ is the intersection of 
all homogeneous systems of parameters consisting of elements of degree $n$.
\end{enumerate}
\end{prop}

\begin{proof}
Suppose that $I$ is generated by elements of degree $n$ and let 
$w$ be a homogeneous element of $\overline I$. Then $w$ satisfies a homogeneous
equation of integral dependence
$$
w^t + a_1 w^{t-1} + \dots + a_{t-1}w + a_t = 0
$$
where $a_i \in I^{i}$. Since $I^i$ is generated by elements
of degree $in$, we see that the degree of $w$ must be at least $n$:
otherwise the homogeneity forces all $a_i = 0 $ and so $w$ would be nilpotent.
Thus no element of degree less than $n$ can be in $\overline I$.

On the other hand,  fix any system of parameters $\{x_1, \dots, x_d\}$ 
consisting of elements
of degree $n$. 
 Then we have a finite integral extension
of graded rings 
$$
A = k[x_1, \dots, x_d] \hookrightarrow S,
$$
and so every homogeneous element $w$ of $S$
satisfies a homogeneous equation of integral dependence
$$
w^t + a_1 w^{t-1} + \dots + a_{t-1}w + a_t = 0,
$$
where $a_i \in A$. 
Now if $w$ has degree at least $n$, the elements $a_i$ have degree 
at least $in$. Thus $a_i \in (x_1, \dots, x_d)^i$, and $w$ is
integral over $(x_1, \dots, x_d)$.

Let $J$ be any minimal 
 homogeneous reduction of $I$ (or of $\overline I$).  From the
above arguments, $J$ is necessarily generated by degree $n$ elements, and
since the ground field is infinite, these generators       form a
system of parameters.
This completes the proof. 
\end{proof}

\subsection{Canonical sheaves and trace.}\label{can}
A general reference for the material in this section is \cite{RD}.
See also \cite{con}.

A canonical sheaf 
for  a Noetherian scheme of dimension $d$ is defined to  be the coherent
sheaf given by the
$-d$-th cohomology of a normalized dualizing complex for the scheme,
when such a dualizing complex exists. 
When the scheme is Cohen-Macaulay, the canonical sheaf is a dualizing
sheaf in the sense of Grothendieck. 
 Dualizing complexes
exist for any equidimensional
 scheme essentially of finite type over an affine Gorenstein 
scheme;
see 
\cite[p299, p306]{RD}.

The canonical sheaf of a scheme $Y$ is not 
uniquely determined up to isomorphism in general. However,
in many situations, there is a canonical choice for the
canonical sheaf. For example, if $Y$ is a normal algebraic variety,
then the usual notion of the canonical sheaf (namely, the unique reflexive
sheaf that agrees with the sheaf of top differential K\"ahler forms on the
smooth locus; see (\ref{nonvan}))
provides a natural, ``truly canonical'' choice for $\omega_Y$. 
More generally, when all schemes are equidimensional and 
 essentially of finite type
over a fixed ground  
 scheme $\spec A$ that admits a residual complex, there is a 
functorial procedure (``upper shriek,'' denoted ``!'')
 for constructing  natural dualizing complexes on them from the
given one on $A$.
In this paper, 
 all schemes will be of essentially finite type over 
a fixed local ground scheme possessing a dualizing complex
 (a field in the geometric settings),
 with respect to which
 all canonical sheaves
will be constructed. 

 For affine schemes, we  also use the terminology
{\it canonical module.}
The canonical module of a local ring $(A, \fm)$
can be described as a finitely generated module whose dual is
the top local cohomology module $H^d_{\fm}(A)$, where $d = \dim A$
(with the duality as described in Subsection 2.3).

\subsubsection{The Canonical sheaf and Cohen-Macaulayness.}
The canonical sheaf of a Cohen-Macaulay scheme is itself a Cohen-Macaulay 
sheaf  of the same dimension. 
Even for non-Cohen-Macaulay schemes, the canonical sheaf
retains a bit of its Cohen-Macaulayness: the canonical sheaf
satisfies Serre's $S_2$ condition. 
 In particular, the canonical sheaf of
a scheme $Y$ of dimension at most two is a Cohen-Macaulay sheaf, even if
$Y$ itself is not Cohen-Macaulay. The
 basic properties of 
 canonical modules over (non-Cohen-Macaulay) local rings are summarized,
with references, in
 \cite[Section 2]{HH}.

\medskip
\subsubsection{Trace.}\label{trace}
Fix a 
proper map
$$Y \overset{f}{\la} \spec A,$$
 where $\spec A$ is an equidimensional local scheme of dimension $d$  with 
fixed  normalized residual complex  $\cR^{\bullet}$.
Letting 
$\cR_Y^{\bullet} := f^! \cR_{\bullet}$ be the corresponding
normalized residual complex for $Y$,
Grothendieck's general theory of {\it trace}
provides a map
$$
\bR\Gamma(Y, \cR_Y^{\bullet}) \la \cR^{\bullet};
$$ 
see  \cite[p318]{RD} for details on $f^!$  and \cite[p383]{RD} 
 for more on 
trace. 
In particular, when $A$ and $Y$ have the same dimension $d$, we can take the
$-d$-th cohomology, and we get a natural 
 map of $A$-modules
\begin{equation*}
\Gamma(Y, \omega_Y) \la \omega_A,
\end{equation*}
which we call the {\it trace map}.
 The dual of the trace map is the natural map of local cohomology modules
$$
H^d_{\fm}(A) \la H^d_Z(\cO_Y),
$$
where $\fm$ is the unique maximal ideal of $A$ and $Z = f^{-1}(\{\fm\})$
(with the duality as discussed below in Subsection \ref{duality}). 

\medskip
The most important case for us is when 
  $Y \la \spec A$ is proper and
birational, specifically, a blowup along an ideal in $\spec A$.
 In this case,
 the trace map turns out to be {\it injective};  dually, 
 the natural map of local cohomology modules $
H^d_{\fm}(A) \la H^d_Z(\cO_Y)
$
is surjective. See, for example, \cite[p103]{LT}.
Throughout this paper, we will always identify 
$\Gamma(Y, \omega_Y) $ with a submodule of $\omega_A$ using this
trace morphism when $Y \la \spec A$ is  projective and birational.

\medskip
In geometric situations, the meaning of the trace map is often clear.
For example,  if both $Y$ and $\spec A$ are normal and $Y \overset{f}{\la} 
\spec A$ is   the blowing up of some closed subscheme $W \subset \spec A$
of codimension at least 2, then there is
a natural inclusion 
$$
\Gamma(Y, \omega_Y) \hookrightarrow \omega_A
$$
obtained by restricting a global section of $\omega_Y$ to $Y \setminus
 f^{-1}(W) = \spec A \setminus W$ and then extending uniquely
to a section of $\omega_A$.  

\medskip
\subsubsection{Adjunction.}\label{adjunction}
Let $E$ be a subscheme of a Cohen-Macaulay 
scheme $Y$  defined locally by a 
nonzerodivisor. Let $\cO_Y(-E)$ denote the ideal sheaf of $E$, so that 
  $\cO_Y(E)$ denote its dual. Then the ``upper shriek'' 
construction in which a canonical module on $Y$
determines one on $E$ is easy to describe. Explicitly,
$$
\omega_E := \SExt^1_{\cO_Y}(\cO_E, \omega_Y).
$$
In particular, $\omega_Y$ and $\omega_E$ are related by the following
exact sequence
$$
 \quad 0 \la \omega_Y \la \omega_Y(E) \la \omega_E \la 0,
$$
 obtained from the long exact sequence that 
arises by applying  the functor $\SHom_{\cO_Y}(-, \omega_Y)$
to the exact sequence $
 0 \la \cO_Y(-E) \la \cO_Y \la \cO_E \la 0.$

If $Y$ fails to be Cohen-Macaulay but is proper over a
local scheme, and $E$ is locally defined by a 
nonzerodivisor, then we still have an exact sequence 
$$
 \quad 0 \la \omega_Y \la \omega_Y(E) \la \omega_E,
$$
but exactness on the right can fail in general. 

\medskip

\subsection{Duality}\label{duality}
Throughout this paper, the term ``duality'' always refers to the following
version of Grothendieck duality combining  global 
and local  duality as developed in \cite{RD}. For a careful proof of this
form of duality, see \cite[p188]{L1}.

\medskip
Let $(A, \fm)$ be a local ring that is a homomorphic image of
a Gorenstein ring (for example, essentially of finite type over a field).
Let $Y \overset{f}{\la} \spec A$ be a proper morphism and let 
$Z = f^{-1}(\{\fm\})$ denote its closed fiber. 
If $Y$ is Cohen-Macaulay of equidimension $d$, then for 
any coherent $\cO_Y$-module $\cF$, there exist $A$-module
isomorphisms for all $i$
\begin{equation}\label{eq3}
H^i_{Z}(Y, \cF) \overset{\cong}{\la} 
\Hom_A(\Ext_Y^{d-i}(\cF, \omega_Y), E_A(A/\fm)),
\end{equation}
where
$E_A(A/\fm)$ is an injective hull of the residue field of $A$,
and $ H^i_{Z}(Y, \cF)$ denotes the local cohomology module of $\cF$ 
with supports in $Z$.  
In particular, if $\cF$ is invertible, then 
$$
H^i_Z(Y, \cF) \qquad {\text{ is dual to }} \qquad
 H^{d-i}(Y, \cF^{-1}\otimes\omega_Y).
$$

Note that this duality includes Serre duality as a special case.
Indeed,  if $Y$ is  a projective variety over a field $k$,
then applied to the proper map $Y \la \spec k$ 
(so that $Y = Z$ and $E_A(A/\fm) = k$), we recover the standard statement
of Serre duality (as in \cite[p243]{Ha}).
At the other extreme, taking $f$ to be the identity map $\spec A \la \spec A$,
we recover the standard local duality familiar to commutative algebraists
(as in \cite[p133]{BH}). 

When $Y$ is not Cohen-Macaulay, the isomorphism (\ref{eq3})
holds as stated only for $i = d$. To get duality for all $i$, 
one must replace $\omega_Y$ in the statement of  (\ref{eq3})
by the normalized dualizing complex for $Y$. We will not
need such general formulations of duality in this paper.
See \cite[p188]{L1}.

\subsection{Rees Rings and Associated Graded Rings.}\label{rees}
Let $I$ be an ideal in a Noetherian ring $A$. 
The {\it Rees ring} of $A$ with respect to $A$ is the $\bbN$-graded ring
$$ A[It] := A \oplus I \oplus I^2 \oplus I^3 \oplus \dots,
$$
and the {\it associated graded ring} or {\it form ring} is the $\bbN$-graded
 ring
$$
 \gr_I A := A/I \oplus I/I^2 \oplus I^2/I^3 \oplus I^3/I^4 \oplus \dots.
$$
In both cases, the ``multiplication'' is the one naturally induced by
multiplication in $A$.
If the ideal $I$ has positive height, then the Rees ring has dimension
$d + 1$, where $d = \dim A$, and the associated graded ring has dimension $d$.
(See {\it e.g.} \cite[Theorem 9.7]{HIO}).

\medskip
Now let $R$ and $G$ be the Rees and associated graded rings, respectively,
 for
$A$ with respect to some ideal $I$ of positive height.
Set $Y = \proj R$. By definition,
the natural projection $Y \overset{\pi}{\la} \spec A$
(induced by the inclusion of $A$ in $R$) is
the blowing up morphism of the ideal $I$. The ideal sheaf $I\cO_Z$ is
invertible, and  defines the scheme theoretic pre-image  of
the subscheme of $\spec A$ defined by $I$.  Thus this pre-image is a
divisor, $E$, called the {\it exceptional divisor}
{\footnote{Caution: If $I$ has height one, the actual exceptional set for the
map $\pi$ may not
be a divisor at all, but a proper subset of $E_{red}$.}}  of $\pi$.
The natural isomorphism 
$R \otimes_A A/I \la G$
identifies $E$ with $\proj G$, so there is a fiber diagram
$$
\begin{CD}
Z @>>> E := \proj G  @>>> Y := \proj R \\
@VVV @VVV   @VVV \\
\spec (A/\fm) @>>> \spec (A/I) @>>> \spec A,
\end{CD} 
$$
where $Z$ is the scheme-theoretic fiber over the closed point $\fm$ of 
$\spec A$. In this diagram, the horizontal maps are all closed
embeddings whereas the vertical maps are all proper.
 If $I$ is an $\fm$-primary ideal,
the schemes $Z$ and $E$ share the same reduced subscheme. 
The  invertible sheaves $I^n\cO_Y$ 
can be identified with the coherent
 sheaves $\cO_Y(n)$ arising from the graded $R$-modules $R(n)$ (where
$R(n)_m = R_{n+m}$).
This justifies our 
 use  of the notations
$$
I^n\cO_Y, \qquad \cO_Y(n), \,\,\, {\text { and }} \,\,\, \cO_Y(-nE)
$$
interchangeably, even when $n $ is negative.

\medskip

\subsubsection{Arbitrary Filtrations.} 
More generally, Rees rings and associated graded rings can be defined
with  respect to an arbitrary filtration of a Noetherian ring $A$.
Let $\{I_n\}_{n \in \bN}$ be a filtration of $A$, that is,
a descending sequence of ideals satisfying $I_n I_m \subset I_{n+m}$ for all
$n, m \in \bN$ and $I_0 = A$.
Then the Rees ring and associated graded rings are defined by
$$
A \oplus I_1 \oplus I_2 \oplus \dots \,\,\,{\text { and }} \,\,\,
 A/I_1 \oplus I_1/I_2 \oplus I_2/I_3 \oplus \dots,
$$
respectively.
The standard Rees and associated graded rings of the previous paragraph
correspond to the filtration $\{I_n\} = \{I^n\}$. In general, the 
Rees ring and associated graded ring of an arbitrary filtration need not
be Noetherian. However, if they are Noetherian and $I_1$ has positive
height, then the dimensions are $d+1$ and $d$, respectively, where $d$ is
the dimension of $A$, just as for filtrations by powers of  ideals.

Given any filtration  $\{ I_n\}$ and a fixed natural number $k$,
there is a Veronese sub-filtration whose $n$-th member is $I_{kn}$.
In this case, the corresponding Rees ring is the $k$-th Veronese  
subring of the Rees ring for $\{ I_n\}$. (The affect on the
associated graded ring is more subtle.) Because
every finitely generated graded algebra has a Veronese sub-ring
generated in degree one,   any filtration giving rise to 
a finitely generated Rees ring  admits a Veronese sub-filtration 
consisting of powers of an ideal $I_k$.

\medskip
The only type of filtration we use in this paper (other than powers
of $I$) is the ``natural'' filtration in a graded
ring: If $S$ is a Noetherian $\bN$-graded ring, then set $I_n = S_{\geq n}$ 
to be the ideal generated by elements of degrees at least $n$.
In this case, the Rees ring, denoted $S^{\natural}$,  is Noetherian
and the associated graded ring is naturally isomorphic to $S$.

\smallskip
\subsubsection{The a-invariant.}\label{a-invariant}
Let 
 $\cR$ be an arbitrary Noetherian
 graded ring over a local ring, and let
$\fM$ 
 denote its unique homogeneous maximal ideal.
  The {\it $a$-invariant} of $\cR$  is defined as
$$
a(\cR) = \max_{n} \{[H^{d}_{\fM}(\cR)]_n \neq 0 \}. 
$$
The a-invariant of a Rees ring is always 
$-1$ while the $a$-invariant of the associated graded ring carries
subtle information about the singularities of $A$ and $R$ 
(see \ref{SSfacts} below).
The term ``a-invariant'' is due to Goto and Watanabe;
 see \cite{GW}.

\medskip

\subsection{The Sancho de Salas Sequence.}\label{SS}
Let 
$$
R = R_0 \oplus R_1 \oplus R_2 \oplus \dots
$$
 be an arbitrary graded ring over a ring $R_0 = A$,
and let $\fm$ be an arbitrary ideal of $A$. 
Set $Y = \proj R$ and $Z = Y \times_{\spec A} \spec A/\fm$.
Then for any graded $R$-module 
 $N = \bigoplus_{n\in \bbZ} N_n$,
 there is a  degree-preserving long exact sequence:
$$
\dots \la H^i_{\fm_R}(N) \la 
\bigoplus_{n \in \bbZ} H^i_{\fm}(N_n) \la 
\bigoplus_{n \in \bbZ} H^i_{Z}(Y, \cN_n) \la
H^{i+1}_{\fm_R}(N) \la \dots
$$
where $\fm_R = \fm \oplus R_1 \oplus R_2 \oplus \dots$
 and $\cN_n$ denotes the quasi-coherent $\cO_Y$-module
corresponding to the graded $R$-module $N(n)$.
This very useful sequence was introduced in \cite{SS} in a special case,
 and later developed by  Lipman in \cite{L2}.

\bigskip
\subsubsection{Local vs Global Cohomology}\label{loc-glob}
For example, consider the extreme case where $\fm = 0$.
Then $Y = Z$ and because $H^i_0(N) = 0$ for $i > 0$, 
the Sancho de Salas sequence degenerates to the
 long exact sequence
$$
0 \la H^0_{R_{>0}}(N) \la 
N \la 
\bigoplus_{n \in \bbZ} H^0(Y, \cN_n) \la
H^{1}_{R_{> 0}}(N) \la 0,
$$
where ${R_{> 0}} = R_{1} \oplus R_2 \oplus \dots $ 
is the ``irrelevant ideal'' of the graded ring $R$, and the 
graded isomorphisms
\begin{equation*}\label{eq}
\bigoplus_{n \in \bbZ} H^i(Y, \cN_n)  \cong H^{i+1}_{R_{>0}}(N)
\qquad{\text{ for }} i \geq 1.
\end{equation*}
This is the familiar identification between sheaf cohomology 
 on a projective scheme and the corresponding local cohomology 
with supports in the irrelevant ideal.
\bigskip

\subsubsection{The case of  Rees rings}\label{SSfacts}
Let 
 $R$ be the  Rees ring 
of a local ring $(A, \fm)$ of dimension $d$ 
with respect to a Noetherian filtration of ideals $I_n$ of 
positive height. In this case, 
$\fm_R = \fm \oplus R_1 \oplus R_2 \oplus \dots$ is the 
unique homogeneous maximal ideal $\fM_R$ of $R$. 
For the case $N = R$, 
 the 
Sancho de Salas sequence is 
$$
\dots \la H^i_{\fM_R}(R) \la  \bigoplus_{n\in \bbN} 
H^i_{\fm}(I_n) \la 
\bigoplus_{n \in \bbZ} H^i_{Z}(Y, \cO_Z(n)) \la
H^{i+1}_{\fM_R}(R) \la \dots.
$$
This exact sequence  can be used to quickly deduce many 
useful well-known
 facts, including:

\begin{enumerate}
\item
As graded $A$-modules, 
$H^{d+1}_{\fm_R}(R) \cong \oplus_{n < 0}H^d_Z(Y, \cO_Y(n))$. This is because 
the maps $H^d_m(I_n) \la H_Z^d(Y, \cO_Y(n))$ are surjective for all
 $n \geq 0 $
(see, e.g.~ \cite[p. 103]{LT}.)
\item
If $R$ is Cohen-Macaulay, then $ H^i_{Z}(Y, \cO_Y(n)) = 0$ for all
$n < 0$ and all $i < d$. 
By duality, this is the same as
$  H^i(Y, \omega_Y(n)) = 0 $ for all $n > 0 $ and all $i >0$.
\item If $I$ is $\fm$-primary and 
both $A$ and $R$ are  Cohen-Macaulay,
 then  also 
$H^i_Z(Y, \cO_Y(n)) = 0 $ for all $n\geq 0 $ and all $1 <i < d$.
  The dual statement is 
$H^i(Y, \omega_Y(n)) = 0$ for all $n\leq 0 $ and all $ 0 < i <  d-1$.
\item If $A$ is Cohen-Macaulay, then $R$ is Cohen-Macaulay if
 and only if $G$ is Cohen-Macaulay with negative $a$-invariant.
 \cite{GN}, \cite{L2}
\item If $G$ is Cohen-Macaulay, then $H^i_Z(Y, \cO_Y) = 0 $ for all $i<d$.
Dually, $H^i(Y, \omega_Y) = 0$ for all $i >0$. 
\cite{SS}, \cite{L2}
\end{enumerate}
To deduce the above statements involving Cohen-Macaulayness, 
use the fact that a module $M$ over a local (or graded)
ring $(\cR, \cM)$ is Cohen-Macaulay if and only if
the local cohomology modules $H^i_{\cM}(M)$ vanish for all $i < \dim M$.

\subsection{The Graded Canonical Module.}\label{gradcan}
Let $R = \bigoplus_{n \in \bbN} R_n$ be an $\bbN$-graded ring finitely 
generated over a local ring $R_0 = A$, where $(A, \fm) $
 is a homomorphic
image of a Gorenstein local ring. 
Let $\fM_R = \fm \oplus R_1 \oplus R_2 \oplus \dots$ denote
the unique homogeneous maximal ideal of $R$. 
Then $R$ admits a {\it graded canonical module,}
which by definition, is a finitely generated graded $R$-module
such that $$
\underline{\Hom}_A(\omega_R, E_A(A/\fm)) \cong H^{\dim R}_{\fM_R}(R),
$$
where 
$E_A(A/\fm)$ is an injective hull of the $A$-module $A/\fm$ and 
the notation $\underline{\Hom}$ denotes
``graded homomorphisms'', namely
 $$\underline{\Hom}_A(\omega_R, E_A(A/\fm)) = 
\bigoplus_{n\in \bbZ}\Hom_A({[\omega_R]}_{-n}, E_A(A/\fm)).$$
In other words,
 a graded canonical module is a finitely generated graded $R$-module
whose ``graded Matlis dual''  is isomorphic, as a graded $R$-module,
 to the top local cohomology 
module with supports in the unique homogeneous maximal ideal of $R$. 
The graded canonical module is uniquely determined up to degree preserving
homomorphism. Furthermore, it is a canonical module for $R$ in the non-graded 
sense as well.  
For details and generalities on graded canonical modules
and related material, see \cite{HIO}, Chapter VII, especially
Section 36, or \cite{BH}, Section 3.6, or the original paper of 
Goto and Watanabe \cite{GW}.

\subsubsection{The canonical module for rings graded over a field.}\label{canS}
Let $S$ be an $\bbN$-graded ring over a field $k = A$,  and let $X = \proj S$
be the corresponding projective scheme. By definition, the graded pieces
of the graded canonical module $\omega_S$ are dual to the graded pieces
of $H^{\dim S}_{\fM_S}(S)$. In the case where $X$ has dimension at least one
(so $S$ has dimension at least two), this latter module can be identified with 
(see \ref{loc-glob})
$$
\bigoplus_{n \in \bbZ} H^{\dim X} (X, \cO_X(n)).
$$
So, 
 we  use duality 
(which holds at the top spot even when $X$ is not Cohen-Macaulay)
to conclude that
$$
\omega_S = \bigoplus_{n\in \bbZ} H^0(X, \omega_Y(n))
$$
is a graded canonical module for $S$.

\subsubsection{The canonical module for $R$.}\label{canR}
 In the  case where $R$ is a Rees ring of an 
 ideal of positive height in a local ring $(A, \fm)$, 
we have 
\begin{equation*}
\omega_R = \bigoplus_{n > 0} H^0(Y, \omega_Y(n)),
\end{equation*}
where 
$\omega_Y$ is the canonical module on $Y = \Proj R$
 constructed from the
fixed one on $A$. Indeed, its dual is 
$
\bigoplus_{n < 0} H^d_Z(Y, \cO_Y(-n)),
$
which is identified with $H^{d+1}_{\fM_R}(R)$ as a graded module, by
\ref{SSfacts} (1). More generally, this argument shows that for 
any graded ring $R$ over a local ring $A$, $[\omega_R]_n = H^0(Y, \omega_Y)n)$
for positive $n$, where $Y = \Proj R$. 

\subsubsection{The canonical module for $G$.}\label{canG}
Likewise, there is a similar choice for the associated graded ring $G$ of
 an $\fm$-primary ideal,
at least when $ d \geq 2$, namely 
\begin{equation*}
\omega_G = \bigoplus_{n \in \bbZ} H^0(E, \omega_E(n)),
\end{equation*}
where 
$\omega_E$ is the canonical module on $E = \proj G$ constructed from the
fixed one on $A$. Indeed, its dual is 
$
\bigoplus_{n \in \bbZ} H^{d-1}(E, \cO_E(-n)),
$
which is identified with $H^{d}_{\fM_G}(G)$ as a graded module, by
the identity (\ref{loc-glob}) above.  In dimension one, 
these arguments give only that 
$$
[\omega_G]_n = 
H^0(E, \omega_E(n))
$$
for {\it positive}  $n$.
The non-positive pieces of $\omega_G$ are more complicated to describe
(and are of crucial importance in our arguments). These 
 will be treated later in 
Lemma \ref{dimonecase}.

\medskip

\subsubsection{The effect of killing a parameter}\label{hyperfact}
Let $\omega_G$ be the graded canonical module of a Cohen-Macaulay graded ring
$G$ over a local ring $G_0$,
 and let $x$ be any homogeneous nonzerodivisor.
Set $\overline{G}$ to be the graded ring $G/xG$. 
Then $$
\omega_{\overline{G}} \cong  \frac{\omega_G}{x\omega_G}(n),
$$ 
as graded modules, where
$n$ is the degree of $x$. 
Here recall that for a graded module $M$, 
 the notation $M(n)$ denotes the same module with the degree shifted so that
$[M(n)]_t = [M]_{n+t}$.
This is well-known and 
easy to  prove; see, for example, \cite[Corollary 3.6.14]{BH}.

\medskip
More generally, even when $G$ is not Cohen-Macaulay, we often 
have useful statements along these lines.
For example, we will make use of the following
 proposition:
\begin{prop}\label{hyper}
Let $G$ be an $\bN$-graded ring of dimension $d>0$, finitely 
generated  over an Artin local ring
$G_0$, and let $x$ be any homogeneous element of degree $n$.
Then there is a natural degree-preserving  injection
$$
\frac{\omega_G}{x\omega_G } (n) \hookrightarrow \omega_{\overline G}, 
$$
where $\overline G$ denotes the ring $G/xG$, 
whenever the dimension of $\ann_G(x)$ (as a $G$-module)
is
 strictly
less $d$.
\end{prop}

\begin{proof}
Consider the four term exact sequence of degree preserving maps
$$
0 \la  \ann_G(x) \la G(-n) \overset{x}{\to} G \la G/xG \la 0. 
$$
Breaking this up, the short exact sequence
$$
0 \la  \ann_G(x) \la G(-n) \overset{x}{\to} xG \la 0
$$
induces an isomorphism
$$
H^d_{\fM_G} (G)(-n) \overset{x}{\la} H^d_{\fM_G}(xG),
$$
since $H^d_{\fM_G}(\Ann_{G}(x)) = 0$.
So the sequence 
$$
0 \la xG \la G \la G/xG \la 0
$$
gives rise to  the exact sequence
$$
H^{d-1}_{\fM_G}(G/xG) \la H^d_{\fM_G}(G(-n)) \overset{x}{\la}
 H^d_{\fM_G}(G) \la 0.
$$
Dualizing  (that is, applying the graded Matlis dual functor; see
 \ref{gradcan})
yields
an exact sequence
$$
0 \la \omega_G \overset{x}{\la} \omega_G(n) \la \omega_{G/xG},
$$
which provides the natural inclusion
$$
\frac{\omega_G}{x\omega_G}  (n) \hookrightarrow \omega_{\overline G}.
$$
\end{proof}

\section{The Main Technical Theorem}

This section forms the technical heart of 
this paper. Here we prove Theorem \ref{mainalg},
 which will later be used to deduce
 the higher dimensional
versions of the Huneke-Swanson formulas for core, our
 ``local version''
of Kawamata's conjecture,
  a formula for 
core conjectured by Corso, Polini and Ulrich,  the results 
linking core and coefficients ideals, and other properties of core in
Section 5 and 6.4.

\smallskip
Recall that 
 if 
$Y = \proj A[It] \la \spec A$ is the blow-up of 
$\spec A$ along an ideal $I$, then the $A$-module
$H^0(Y, \omega_Y)$ can be naturally identified with a 
submodule of $\omega_A$ (see \ref{trace}).
Likewise,  the modules $H^0(Y, I^n\omega_Y)$ for $n \geq 0$ can 
 be identified with submodules of $\omega_A$, denoted by $\Omega_n$.

\begin{thm}\label{mainalg} 
Let $I$ be an $\fm$-primary ideal in 
a local ring $(A, \fm)$ of positive dimension $d$ 
 containing the rational numbers.
 Assume that for any reduction $J$ of $I$,
\begin{equation}\label{hyp} 
J\omega_A \cap \Omega_{d-1}  = J (\Omega_{d-2} \cap \omega_A),
\end{equation}
as submodules of $\omega_A$.
Then 
\begin{equation*}
\core{(I\omega_A)} \subset \Omega_d 
\end{equation*}
as submodules of $\omega_A$, where 
$\core{I\omega_A}$ denotes the intersection in $\omega_A$
of the submodules $J\omega_A$ as $J$ ranges through  
 all reductions of $I$. 
\end{thm}

\begin{rem}\label{one}
Note that if $d \geq 2$, then $\Omega_{d-2}$ is contained
in  $\omega_A$ in any case; the intersection with $\omega_A$ is
relevant only when $d =1$.
\end{rem}

\begin{rem}\label{304}
For ideals $I$ of reduction number at most one,
 the assumption that $A$ contains the rational
numbers is unnecessary; see Remark \ref{d=2}.

\end{rem}
\begin{rem} In the geometric setting, the module $ H^0(Y, I^d\omega_Y)$
is  closely related to the adjoint ideal (or multiplier ideal)
of $I^d$, at least when 
 $Y$ happens to be smooth (or have rational singularities).
See Remark \ref{mult.ideal}.
\end{rem}

 As we will see in the next section, the hypothesis (\ref{hyp}) of Theorem
\ref{mainalg}
is a type of  ``Brian\c con-Skoda'' statement, and  it
is satisfied in many nice situations.
For example, we will show that $(\ref{hyp})$ holds  whenever
 $A$ and the Rees ring $A[It]$
 are Cohen-Macaulay, or more generally even if $A$ is
not Cohen-Macaulay provided that the irrelevant ideal of the Rees ring $A[It]$
is a Cohen-Macaulay module. 
This latter condition arises naturally in the geometric setting that
motivates us.

\bigskip
The proof of Theorem \ref{mainalg} will occupy this entire
section. 
The main point is that 
Lemma \ref{key} 
reduces us
to a related statement about the intersection of
the corresponding 
submodules of the canonical module $\omega_G$ of the associated
graded ring $G$. 
This statement about $\omega_G$ is then proved by induction on the 
dimension,   
with the hard part being the case where $d = 1$.
 For all these steps, 
we need a rather delicate understanding of the structure of the
canonical module for $G$---especially how its graded pieces 
are related to the adjoint-type modules  $\Omega_n$. Thus we begin in 
Subsection 3.1 with a detailed study of the modules $\Omega_n$.

\medskip\subsection{The filtration by Adjoint-type modules.}
We fix some notation to be used throughout the  rest of this section.
We let  
 $(A, \fm)$
 denote  a local ring of dimension $d\geq 1$ which is 
assumed to have an infinite residue field 
(and as always, possesses a canonical module).
 Let $I$ denote a
 proper ideal in $A$ of positive height.
 Let $R$ and $G$ denote
the Rees ring and the associated graded ring with respect to $I$, respectively.
We set $Y = \proj R$ and let $ Y \overset{\pi}{\to} \spec A$ denote the natural
blowing up morphism. As always, 
 we  identify 
$\pi_*\omega_Y$ with a submodule of $\omega_A$ (see Subsection \ref{trace}).
\smallskip

We first establish some elementary properties of the 
``adjoint-type'' $A$-modules 
$$
\Omega_n := \Gamma(Y, I^n\omega_Y),
$$
where $n \in  \bbZ$. First note that:
\begin{enumerate}
\item Each $\Omega_n$ is a finitely generated $A$-module.
\item If $n \geq m$, then $\Omega_n \subset \Omega_m$.
\item For all $n \in \bbZ$, we have $ I\Omega_n \subset \Omega_{n+1}.$
\item There are natural identifications
$\Omega_n = \Hom_A(I^p, \Omega_{n+p}) $ for all $ n \in \bbZ$ and all 
$p \geq 0$. 
\end{enumerate}
The first property is immediate from the properness of $Y \la \spec A$,
while the next two properties follow immediately from the definition.
The fourth property follows from the useful but elementary general fact: 
{\it {Let $\cF$ be a coherent $\cO_Y$-module such that the local generator
for $I\cO_Y$ is a nonzerodivisor on $\cF$ at each point of $Y$.
Then $\Hom_Y(J\cO_Y, \cF) = \Hom_A(J, \Gamma(Y, \cF))$ for any ideal 
$J \subset A$.}} (This fact is easy to prove;
see, for example, Lemma 2.1 of \cite{Hy2}.)

\bigskip
Our proof of Theorem \ref{mainalg} will exploit the
following 
 relationship between the $\Omega_{n}$ 
 and the graded pieces of $\omega_G$.

\begin{lem}[Cf. ~\cite{Hy2}, Theorem 2.2e]\label{G}
Let $I$ be an $\fm$-primary ideal in a 
 local ring  $(A, \fm)$ of positive dimension $d.$
Then there is a natural inclusion 
 $$
\bigoplus_{n \geq 1}\, \Omega_{n-1}/\Omega_n 
 \hookrightarrow \omega_G 
$$
of  graded $G$-modules. 
This inclusion is an isomorphism if $Y$ is Cohen-Macaulay and
$H^1(Y, \omega_Y(n) ) = 0 $ for all $n$, for example, if
  both $R$ and $A$ are Cohen-Macaulay. 
\end{lem}

\begin{proof}
Consider the adjunction sequence
$$
0 \la \omega_Y \la \omega_Y(-1) \la \omega_E,
$$
which is exact also on the right if $Y$ is Cohen-Macaulay.
Tensoring with $\cO_Y(n)$ and computing cohomology we get an 
exact sequence of cohomology 
$$
0 \la H^0(Y, \omega_Y(n)) \la H^0(Y, \omega_Y(n-1)) 
\la H^0(E, \omega_E(n)).
$$
Thus 
there is a natural injection 
$$
\bigoplus_{n\in \bbZ} \Omega_{n-1}/\Omega_{n} \hookrightarrow
 \bigoplus_{n\in \bbZ} 
H^0(E, \omega_E(n)) 
$$
for all $n \in \bbZ$.
Since  $\omega_G$ and $ \bigoplus_{n\in \bbZ} 
H^0(E, \omega_E(n))$ agree in positive degrees, 
 (see \ref{canG}), the first claim is proven.

If $Y$ is Cohen-Macaulay, then 
we have an exact sequence
$$
0 \la H^0(Y, \omega_Y(n)) \la H^0(Y, \omega_Y(n-1)) 
 \la H^0(E, \omega_E(n)) \la H^1(Y, \omega_Y(n)), 
$$
so 
the inclusion is a bijection in degree $n$ if 
 $H^1(Y, \omega_Y(n)) = 0$.
When $R$ and $A$ are both Cohen-Macaulay, the scheme $Y = \proj R$ is 
Cohen-Macaulay, and this 
  vanishing holds (see \ref{SSfacts}). The lemma is proved.
\end{proof}

\begin{rem}\label{d>1}
 The proof of Lemma \ref{G} shows that if $d \geq 2$, the inclusion
$
\bigoplus_{n \in \mathbb Z}\, \Omega_{n-1}/\Omega_n 
 \hookrightarrow \omega_G 
$
 holds for all $n$, not just positive $n$.  
A nice  consequence is
 that, for dimension $\geq 2$, 
 the increasing chain of
modules
$$
\dots \subset\Omega_2  \subset \Omega_1 \subset \,   \Omega_0 
\subset \Omega_{-1} \subset  \Omega_{-2} \dots  
$$
must stabilize for $-n < -a$, where $a$ is the
 $a$-invariant of $G$. In fact, it 
 stabilizes to $\omega_A$:
\end{rem}
\begin{lem}\label{stab}
Let $I$ be an $\fm$-primary ideal in a local ring $(A, \fm)$ of
dimension at least two. Then 
 $$\Omega_{-n} = \omega_A$$ for all 
$n > a(G)$. 
\end{lem}

\begin{proof} Fix $n \gg 0$. 
To verify that 
$\omega_A \subset \Omega_{-n},$  
  recall that 
$\Omega_{-n} = \Hom_{A}(I^n, \Omega_0)$. So it will
 suffice to  show that for $n \gg 0$,
 $I^n\omega_A \subset \Omega_0$ (then   
each   $f$ in $\omega_A$ determines
the  element ``multiplication by $f$'' 
in  $\Omega_{-n} = \Hom_{A}(I^n, \Omega_0).$)
Because the blowup map restricts to an isomorphism
away from 
the closed set $\spec (A/I)$,
the trace map  $\Omega_0 \subset \omega_A$ 
becomes
an identity after localizing at any element of $I$. 
So some power of $I$ annihilates $\omega_A/\Omega_0$, and  
$I^n\omega_A \subset \Omega_0$ for large $n$, as needed. 

To check the reverse inclusion, 
consider the 
exact sequence of $A$-modules
$$
0 \la \omega_A \la \Omega_{-n} \la Q \la 0
$$
where, because $I$ is $\fm$-primary, the module
 $Q$ is supported at $\fm$. There is an induced sequence
of local cohomology
$$
 H^0_{\fm}(\Omega_{-n}) \la H^0_{\fm}(Q) \la H^1_{\fm}(\omega_A).
$$
The right term above vanishes because $\omega_A$ is $S_2$ while the
left term vanishes because $\Omega_{-n} = \Hom_A(\Omega_n, \Omega_0)
\subset \Hom_A(\Omega_n, \omega_A)$, so no element can be  killed 
by any parameter
in $A$. Thus $ H^0_{\fm}(Q) = 0,$ and since $Q$ is supported
only at $\fm$, $Q = 0$ as well. 
\end{proof}

In dimension one, the picture is somewhat different:
 although it is still true that $\omega_A \subset \Omega_{-n}$
for large $n$, we
  do not 
get stabilization.

\begin{lem}\label{lem1} If $I$ is integral over a principle ideal, then
 $\Omega_n = x\Omega_{n-1}$ for 
all $n \in \bbZ$, where $x$ generates a minimal reduction for $I$. 
This holds in particular when $I$ is an $\fm$-primary ideal in a local
 ring of dimension one. 
\end{lem}

\begin{proof}
In this case, $Y = \proj A[It]$ is affine, so 
$$\Omega_n = I^n\omega_Y =  x^n \omega_Y = x (I^{n-1}\omega_Y) =
 x\Omega_{n-1}$$
for all $n$.
\end{proof}

\bigskip
The next lemma refines our understanding of the modules
$\Omega_n$ for the case of non-negative  $n$.

\begin{lem}
\label{colonlem}
Let $(A,m)$ be a  local ring and let $I$ be a proper
 ideal of $A$
of  height greater than one. Then $\Omega_{n+1}:_{\omega_A}I=\Omega_n$ 
for all $n\ge 0$.
In fact,
$\Omega_{n+p}:_{\omega_A} I^p =
 \Omega_n $ for all $n \geq 0$ and all $p \geq 1$. 
This also holds for ideals of height one
that are integral over a principle ideal.
\end{lem}

\begin{proof} 
Recall that $\Omega_n = \Hom_{A}(I^p, \Omega_{n+p})$.
Clearly each element $w$ of 
$\Omega_{n+p}:_{\omega_A}I^p$ gives rise to an element of $\Omega_n$,
namely, the ``multiplication by $w$ map'' in  $\Hom_{A}(I^p, \Omega_{n+p})$.
We need to show that every $A$-module map from $I^p$ to $\Omega_{n+p}$
arises this way.

We treat the height greater than one case first.
Because $I^p$ contains a regular sequence of length two on $\omega_A$,
 we have  $\Hom_A(I^p,\omega_A)
=\omega_A$. Since $\Hom_A(I^p,\Omega_{n+p})\subset \Hom_A(I^p,\omega_A)$,
 this implies
that every element of $\Hom_A(I^p,\Omega_{n+p})$ arises by multiplication by 
some element
of $\omega_A$, as needed.

Now suppose 
 $I$ has a reduction generated by one element, say $x$. 
By Lemma \ref{lem1}, 
we have 
$\Omega_{n+1}=x\Omega_n$ for all $n \in \bbZ$.
So also $\Omega_{n+p}=x^p\Omega_n$ for all $n\geq 0$ and all $p \geq 1$.
 Take any $u\in \Omega_n
= \Hom_A(I^p,\Omega_{n+p})$. 
Then $u(x^p)=x^p\omega$ for some $\omega\in \Omega_n \subset \omega_A$.
We claim now that $u$ is the map ``multiplication by $\omega$.''
To check this, take  any $a\in I^p$. We
 get
$x^pu(a)=au(x^p)=ax^p\omega$, so because $x^p$ is a nonzerodivisor
on $\omega_A$, we conclude that $u(a)=a\omega$.
This shows that $\Omega_n = \Hom_A(I^p,\Omega_{n+p}) 
\subset \Omega_{n+p}:_{\omega_A} I^p $.
 The proof is complete.
\end{proof}

\bigskip
\subsection{The canonical module $\omega_G$ in the dimension one case.}
\label{grad-adj}
Because our strategy for  proving Theorem \ref{mainalg} 
is to reduce to the one-dimensional case, 
we need to 
understand $\omega_G$ in the case $d =1$. In dimension one, 
if $R$ and $G$ are both Cohen-Macaulay,
it is easy to see that 
$$\omega_G \cong \bigoplus_{n > 0 } \frac{\Omega_{n-1}}{\Omega_n}.
$$
 However, our reduction to the dimension one case will destroy the
Cohen-Macaulayness of $R$ (even when  it is  assumed at the start),
so the following 
lemma is  crucial.

\begin{lem}\label{dimonecase}
\label{Glem}
Let $(A,\fm)$ be a one-dimensional local ring. Let
$I \subset A$ be an $\fm$-primary ideal. Then 
\begin{equation*}
[\omega_G]_n =
\begin{cases}  \Omega_{n-1}/\Omega_n &  {\text{ for }} n >0 \\
\frac{\Omega_{n-1} \bigcap \Hom_A(I^{-n},\omega_A)}{\Omega_n}
&\text{ for }  n \leq 0.
\end{cases}
\end{equation*}
Here, for $n \leq 0$, 
the intersection is
carried out in  $\Hom_A(I^{-n+1}, \omega_A)$ which 
naturally contains both $\Omega_{n-1} =  \Hom_A(I^{-n+1}, \Omega_0)$
and 
 $\Hom_A(I^{-n}, \omega_A)$. 
\end{lem} 

To understand the 
lemma,
note that for $n \leq 0$, 
the module $\Omega_n$ can be  considered as a submodule of 
$\Hom_A(I^{-n}, \omega_A)$ via the injection 
$\Omega_n = \Hom_A(I^{-n}, \Omega_0) \la  \Hom_A(I^{-n}, \omega_A)$
induced by the inclusion $\Omega_0 \subset \omega_A$.
Also, there is a natural inclusion of $\Hom_A(I^{-n}, {\omega_A})$ in 
$\Hom_A(I^{-n+1}, {\omega_A})$, induced
as the dual of the inclusion  $I^{-n+1} \subset I^{-n}$,
since $\Hom_A(I^{-n}/I^{-n+1}, \omega_A) = 0$.

\begin{proof} 
Consider the long exact sequence 
(see \ref{loc-glob}) relating cohomology for $G$ and $\proj G$:
$$
0 \la
H^0_{\fM_G}(G) \la G \la 
\bigoplus_{n \in \bbZ} H^0(E, \cO_E(n)) \la H^1_{\fM_G}(G)
\la 0.
$$
 Applying 
duality for the proper  map $E \la \spec A/I$, we get an exact sequence
$$
0 \la \omega_G \la \bigoplus_{n\in \bbZ}H^0(E, \omega_E(n)) \la 
\bigoplus_{n\in \bbZ} \Hom_{A/I}(G_n, E_{A/I}) = E_{G},
$$ 
where 
$  E_{A/I}$ denotes the injective hull of the residue field of the zero 
dimensional ring $A/I$ (and $E_G$ is by definition
the graded injective hull of the residue
field of $G$; see \cite[p293]{HIO}). For $n >0$, $[E_G]_n = 0$, so 
there are natural isomorphisms 
$$
[\omega_G]_n =  H^0(E, \omega_E(n))
$$
for all $n > 0$, whereas for $-n \leq 0$, we
have
$$
[\omega_G]_{-n} =
 \ker{[H^0(E, \omega_E(-n)) \la \Hom_A(I^n/I^{n+1}, E_{A/I})]}.
$$ 

Note that $Y$ is always Cohen-Macaulay:
Away from $E$, $Y$ is isomorphic to the zero-dimensional scheme
 $\spec A \setminus \fm$, whereas since $E$ is defined by a nonzerodivisor
locally in $Y$, the Cohen-Macaulay property of the zero-dimensional scheme
$E$ lifts to $Y$. 
 So
we have a twisted adjunction sequence 
 \begin{equation*}0\lrarrow \omega_Y(n) 
\lrarrow \omega_Y(n-1) \lrarrow \omega_E(n) \lrarrow 0.\end{equation*}
Because $Y$ is affine, we have a corresponding sequence
of global sections
\begin{equation*}0\lrarrow \Omega_n \lrarrow 
\Omega_{n-1} \lrarrow 
H^0(E,\omega_E(n)) \lrarrow 0.
\end{equation*}
Hence $\Omega_{n-1}/\Omega_n=H^0(E,\omega_E(n))$ for all $n\in \bbZ$.
So for $n > 0$, we conclude that 
$$
[\omega_G]_n = \Omega_{n-1}/\Omega_n,
$$
as needed.

 It remains to treat the case $n \leq 0$. 
Because the case $n = 0$  is the only one we actually need later, 
we write down the argument carefully only  in this case
for the sake of clarity. However,
the same exact argument ``twisted'' by $\cO(n)$ produces the
result for any $n \leq 0$.

Consider the case $n = 0$. Because 
$\omega_{A/I}$ is an injective hull of the residue field for 
the zero-dimensional local ring $A/I$, 
we have that 
$$
[\omega_G]_0 = 
\ker{[H^0(E, \omega_E) \overset{\eta}{\to} \Hom_A(A/I, \omega_{A/I}) =
 \omega_{A/I}]}.
$$
We need to understand the map $\eta$.
Being dual to the natural map 
$$ 
A/I = H^0_{\fm}(A/I) \la H^0_E(\cO_E)  = \cO_E,
$$
the map $\eta$ is in fact the trace map for the proper map of zero dimensional
schemes $E \la \spec A/I$. (Note that the functors $\Hom_A(-, \omega_{A/I})$
and $\Hom_E(-, \omega_E)$ are identical on $\cO_E$-modules
by the adjointness of tensor and Hom.)

To understand $\eta$, we consider the commutative diagram of $A$-modules,
whose existence we will justify momentarily:
$$
\begin{CD}
0 @>>> \Hom_A(A, \omega_A) @>>> \Hom_A(I, \omega_A) @>>>
\Ext^1_A(A/I, \omega_A) = \omega_{A/I} \la 0 \\
@. @AAA @AAA @A{\eta}AA  @.\\
0 @>>> \Hom_{\cO_Y}(\cO_Y, \omega_Y) @>>> \Hom_{\cO_Y}(I\cO_Y, \omega_Y)
 @>>> \Ext^1_{\cO_Y}(\cO_E, \omega_Y) = \omega_E @>>> 0. 
\end{CD}
$$
In this diagram, the first two upward arrows are injective (note the
first one is the trace map for $Y \la \spec A$), and it is the kernel of
the rightmost upward  arrow that we want to understand.
Knowing that $[\omega_G]_0 = \ker \eta$,
a look at this diagram  shows that  there is a natural identification
$$
[\omega_G]_0 = \frac{\Hom_{\cO_Y}(I\cO_Y, \omega_Y) \cap \Hom_A(A, \omega_A)}
{\Hom_{\cO_Y}(\cO_Y, \omega_Y)} = 
\frac{\Omega_{-1} \cap \omega_A}{\Omega_0},
$$
as claimed. This will complete the proof in the case 
where  $n=0$; the proof for arbitrary negative $n$
is essentially the same, just ``twisted'' by $n$. It remains only to 
justify this commutative diagram.

\smallskip
Finally, 
to justify the diagram is easy.{\footnote{One can do this abstractly, 
using the point of view of \cite{RD}, but we prefer a hands on verification.}}
It is induced from the  diagram of long exact sequences 
 arising
from the natural diagram of $A$-modules
$$
\begin{CD}
0 @>>> I @>>> A @>>> A/I @>>> 0 \\
@. @VVV @VVV @VVV @.\\
0 @>>> I\cO_Y @>>> \cO_Y @>>> \cO_E @>>> 0,
\end{CD}
$$
where the downward arrows are the natural inclusions.
(Recall that $Y$ and $E$ are affine, so we abuse notation,
identifying  $\cO_Y$ and $\cO_E$  with their corresponding rings of global 
sections.)
The existence of a natural induced map
 between the corresponding long exact 
Ext-sequences is a very general fact holding for any diagram of short
exact sequences over any  commutative rings $R \la S$. Indeed, 
say $M$ and $P$ are $R$-modules, and $N$ and $Q$ are $S$-modules.
Given $R$-module maps $M \la N$ and $Q \la P$, there are
naturally induced  functorial maps
$$
\Ext_S^{i}(N, Q) \la \Ext_R^{i}(M, P)
$$
for each $i \in \bbN$. These maps can be viewed as the
composition of three natural maps of $R$-modules:
$$
\Ext_S^{i}(N, Q) \la \Ext_R^{i}(N, Q) \la
\Ext_R^{i}(N, P) \la \Ext_R^{i}(M, P).
$$
The first arrow above is naturally induced by the ``forgetful functor.''
(An injective resolution $Q \la I^{\bullet}$
of $Q$ as an $S$-module can be viewed as an exact sequence of $R$-modules;
then for any injective resolution $Q \la J^{\bullet}$ of $Q$ as an $R$-module,
there will be an induced $R$-module
 map of complexes $I^{\bullet} \la J^{\bullet}$.
This induces an $R$-module map of complexes 
$\Hom_R(N, I^{\bullet}) \la \Hom_R(N, J^{\bullet})$, which can be pre-composed
 with the  map that forgets the $S$-structure
$ \Hom_S(N, I^{\bullet}) \la \Hom_R(N, I^{\bullet})$. This map of complexes
induces a unique map on the level of cohomology, and this is the natural
map of Ext groups we have in the first arrow above.) The second
arrow is the natural covariantness of the functor Ext in the second 
argument, and the third arrow is the natural contravariantness in the first 
argument.

The proof of Lemma \ref{dimonecase} is complete, at least for $n \geq 0$.
 The proof for arbitrary negative $n$
is essentially the same as the proof for $n = 0$, just ``twisted'' by $n$.
(We will use this lemma only in the cases $n = 0$ and $n =1$.)
\end{proof}

\subsection{Reduction to the Associated Graded Ring}
We now state and prove the key lemma  
which 
 provides the crucial step of reducing 
Theorem \ref{mainalg} to a related statement about the
canonical module of the associated graded ring. 

\begin{lem}[Key Lemma]\label{key}
Let $I$ be an $\fm$-primary ideal in a 
local ring  $(A, \fm)$  of  positive dimension $d$.
Let $\cR$ be any non-empty set of minimal reductions of $I$.
As in Theorem \ref{mainalg}, assume that 
$$
J\omega_A \cap \Omega_{d-1} = J(\Omega_{d-2} \cap \omega_A)
$$
 for all reductions $J \in \cR$. 
Then 
$$
\bigcap_{\{x_1, \dots, x_d\} \in \cR} (x_1, \dots, x_d)\omega_A \subset
\Omega_d
$$
if
$$
\bigcap_{\{x_1, \dots, x_d\} \in \cR} 
\left[(x_1^*, \dots, x_d^*)\omega_G\right]_d = 0,
$$
where $x_i^*$ denotes the degree one element of $G$ 
given by the class of $x_i$ modulo $I^2$.
The converse also holds if $Y$ is Cohen-Macaulay and both 
$H^1(Y, \omega_Y(d)) $ and $H^1(Y, \omega_Y(d-1)) $ vanish.
\end{lem}

\begin{proof}
 Set 
$$
W = \bigcap_{\{x_1, \dots, x_d\} \in \cR} (x_1, \dots, x_d)\omega_A.
$$
We already know that $\Omega_d \subset W$, because of our assumption
that $\Omega_d \subset J\omega_A$ for all reductions $J$.
Let $x \in W \setminus \Omega_d$.
Since the only associated prime of $\omega_A/\Omega_d$ is $\fm$,  
we may assume that $x \in \Omega_d:_{\omega_A} \fm
 \subset \Omega_d:_{\omega_A} I.$
But then $x \in \Omega_{d-1}$ by Lemma \ref{colonlem}, and its class in 
$\Omega_{d-1}/\Omega_d $ determines a non-zero element of
$ [\omega_G]_d$ under the natural inclusion $\Omega_{d-1}/\Omega_d 
\hookrightarrow [\omega_G]_d$ guaranteed by Lemma \ref{G}.
 Because 
$J\omega_A \cap \Omega_{d-1} = J(\Omega_{d-2} \cap \omega_A)$
  for any minimal
reduction $J$ of $I$ 
(by hypothesis), we see that 
$$
x \in 
\bigcap_{\{x_1, \dots, x_d\} \in \cR} (x_1, \dots, x_d)(\Omega_{d-2} 
\cap \omega_A).
$$
By Remark \ref{d>1} and Lemma \ref{dimonecase}, we
conclude that 
 the class of $x$ is in 
$$
\bigcap_{\{x_1, \dots, x_d\} \in \cR}[(x_1^*, \dots, x_d^*)\omega_G]_d,
$$
and hence is zero by assumption.
 So in light of the inclusion
$\Omega_{d-1}/\Omega_d \subset \omega_G$, the element $x$ must have
 been in $\Omega_d$ 
after all. The proof that the second condition implies the first is complete.

For the converse, consider a degree 
$d$ element in 
$$
\bigcap_{\{x_1, \dots, x_d\} \in \cR} (x_1^*, \dots, x_d^*)\omega_G.
$$
Our additional hypotheses imply that 
  ${[\omega_G]}_d \cong \Omega_{d-1}/\Omega_{d}$,  so 
this element is represented by some $x$ in $\Omega_{d-1}$.
Likewise, since  ${[\omega_G]}_{d-1} \cong 
(\Omega_{d-2} \cap \omega_A)/\Omega_{d-1}$
and $\Omega_d \subset J\omega_A$, we can assume with out loss of generality
 that
$$x \in (x_1, \dots, x_d)(\Omega_{d-2} \cap \omega_A),$$
for any $(x_1, \dots, x_d) \in \cR$.
Since such $x$ in is 
$
\bigcap_{\{x_1, \dots, x_d\} \in \cR} (x_1, \dots, x_d)\omega_A, 
$
our hypothesis ensures that it is in $\Omega_d$, so that $x$ represents
the zero class in $\omega_G$. This implies that
$$
\bigcap_{\{x_1, \dots, x_d\} \in \cR}[(x_1^*, \dots, x_d^*)\omega_G]_d = 0.
$$
\end{proof}

\subsection{The Proof for $\omega_G$}
In light of the 
 Key Lemma \ref{key},
the next result will complete the proof of 
 Theorem \ref{mainalg}.

\begin{thm}\label{hard}
Let $(A, \fm)$ be a local ring 
 containing the rational numbers.
Let $I$ be any $\fm$-primary ideal. 
Then
$$
\bigcap
_{(x_1, \dots, x_d) \in \cR} [(x_1^*, \dots, x_d^*)\omega_G]_d = 0,
$$
where $\cR$ is the set of all minimal reductions of $I$. 
Furthermore, if $A$ is the localization of a finitely generated $\bbN$-graded
domain $S$ and $I$ is the expansion to $A$  of a homogeneous
 ideal of $S$ generated by elements of all the  same
degree $n$, then the set $\cR$ can be taken to be
the set of all reductions of $I$ generated by elements of degree $n$.
\end{thm}

We will prove Theorem \ref{hard} by induction on $d$. 
The hard part will be to deal with the case where $d=1$.
The inductive step looks slightly technical also, 
because of the necessity of dealing with  non-Cohen-Macaulay associated
graded rings. We isolate most of these technicalities in the following lemma,
which compares the canonical modules of two graded rings
closely related to $G$.  The reader is advised to think about the case where
$G$ is Cohen-Macaulay on a first read through, as this simplifies the
arguments and is sufficient for our main algebraic results (but not
for Kawamata's Conjecture).

\begin{lem}\label{nonCM}
Let $(A, \fm)$ be a local ring 
of dimension $d \geq 2$, and let $G$ be the associated graded ring of $A$
with respect to an $\fm$-primary ideal $I$. Let $y \in I \setminus I^2$
be a general element of $I$, and let $y^* = y + I^2$ 
 denote the corresponding
degree one element  in $G$.
Then 
\begin{enumerate}
\item $(y) \cap I^n  = y I^{n-1}$ for all $n \gg 0$.
\item There is a natural 
 degree preserving surjection
$$
G/y^*G \la \overline G \, := \, 
\frac{\overline A}{\overline I}
\oplus\frac{\overline I}{{\overline I}^2}
\oplus\frac{{\overline I}^2}{{\overline I}^3} \oplus \dots, 
$$
where $\overline G$ denotes the associated graded ring of 
$\overline A = A/yA$ with respect to the image ideal
 $\overline I = I \overline A$, which becomes an bijection in 
degrees $n \gg 0$.   
\item This surjection induces a 
 degree preserving  isomorphism
$$
\omega_{\overline G} \la \omega_{G/y^*G}.
$$ 
\end{enumerate}
\end{lem}

\begin{rem}
When the associated graded ring $G$ is Cohen-Macaulay (as it is in 
our main applications), 
the equality holds in (1) for all $n > 0$ and the map
in (2) is an isomorphism, making (3) obvious. 
\end{rem}

\begin{proof}
Condition (1) is easy to prove and 
well-known: the point is to choose $y$ 
so that $y^*$ avoids the {\it relevant} associated primes of
$G$. (Such an element $y$ is called ``{\it filter regular}.'')
 See, for example, \cite[Lemma 3.2]{Trung}.

Condition (2) follows immediately. 
Indeed, one simply verifies that in degree $n$, this map looks like
$$
\frac{I^n}{xI^{n-1} + I^{n+1}} \la \frac{I^n}{(x)\cap I^{n} + I^{n+1}}.
$$
This is obviously surjective for all $n$, and becomes an
isomorphism  when $n \gg 0.$

To prove (3),  note that the kernel of the 
 natural degree preserving surjection
$$
G/y^*G \la  \overline G
$$
is non-zero in only finitely many degrees,
and  therefore  
has finite length. In particular, 
 the kernel  $K $ is
a zero dimensional graded $G/y^*G$-module. As both graded rings above
have dimension $(d-1)$, the corresponding long exact sequence of
local cohomology induces an isomorphism
$$
H^{d-1}_{\fM_{G/y^*G}}(G/y^*G) \la 
H^{d-1}_{\fM_{\overline{G}}}(\overline G),
$$
since $H^{d-1}_{\fM_G}(K) = 0$.

Dually 
(after applying the graded Matlis dual; see Subsection \ref{gradcan}),  
 we have the
desired degree preserving isomorphism
$$
\omega_{ \overline{G}} \la \omega_{ G/y^*G}.
$$
\end{proof}

\medskip
\begin{proof}[Proof of Theorem \ref{hard}]
 We first carry out the inductive step.
Assume that $d > 1$ and that the theorem has been proven for rings $A$ 
of dimension $d-1$. Suppose that  
\begin{equation*}\omega\in 
\bigcap_{(x_1,\ldots,x_d)\in \cR}[(x_1^*,\ldots,x_d^*)\omega_G]_d.
\end{equation*}
We will construct
a sequence $y_1,y_2\, \ldots$ of elements of $I$ (of degree $n$ in the
graded case) such that the 
corresponding elements $y_1^*,y_2^*\, \ldots$ of $G_1$ are nonzerodivisors
on $\omega_G$
 with the property 
that 
$$
\omega \in (\prod_{i=1}^t y_i^*)\omega_G
$$
 for all $t\ge 1$. This will imply the claim, because 
then $\omega \in \omega_G$ would be divisible by elements of arbitrarily
large degree, which is impossible since $\omega_G$ is Noetherian,
and hence vanishes in sufficiently
negative degrees.

To construct the sequence of $y_i$'s, we proceed inductively.
 Assuming that 
$y_1,\ldots,y_{t-1}$ have already been constructed,
choose a general  $y_t\in I$ (of degree $n$ in the graded case)
 such that $\{y_1^*\cdots  y_{t-1}^*, y_t^*\}$
 is part of a  sequence of parameters for
 $G$. This is possible, since the elements $x^*$ with
$x\in I$ (of degree $n$) generate an $\fM_G$-primary ideal. 
We must verify that $\omega$ is divisible by the product 
$y_1^*\cdots  y_{t-1}^*y_t^*.$

\medskip
Consider the class of 
 $\omega$ modulo $y_t^*\omega_G$. 
This is an element of $\omega_G/y_t^*\omega_G$ of degree $d$,
so in light of the injection provided by Proposition \ref{hyper},
it determines an element of degree $d-1$ in 
$\omega_{\frac{G}{y_t^*G}}$. By Lemma \ref{nonCM}, 
therefore, we can interpret it as an element
in 
$$
[\omega_{\overline G}]_{d-1},
$$
 where $\overline G$ is the associated graded ring of the ring 
$A/y_tA$ with respect to the  image of $I$.

Now any minimal reduction of $\overline{I}$ 
in  $\overline A$, say 
 $(\overline{x_1},\ldots,\overline{x_{d-1}})$, lifts to a minimal
reduction
$(x_1,\ldots,x_{d-1}, y_t)$ 
of $I$,
because, again by Lemma \ref{nonCM}, for all $r \gg 0$, 
 \begin{equation*}I^{r+1}\subset 
(x_1,\ldots,x_{d-1})I^r+(y_t)\cap I^{r+1}\subset (x_1,\ldots,x_{d-1})I^r+
y_tI^r.
\end{equation*}
It thus follows that
\begin{equation*}
\overline{\omega} \in \bigcap_{\overline{x_1},\ldots,\overline{x_{d-1}}}
[(\overline{x_1}^*,\ldots,\overline{x_{d-1}}^*)\omega_{\overline G}]_{d-1}=0,
\end{equation*} 
where ${\overline{x_1},\ldots,\overline{x_{d-1}}}$ ranges through the 
set of all minimal reductions of $\overline{I}$ in $\overline A$.
(In the graded case, we assume the $x_i$'s all to have degree $n$.)

By the inductive 
hypothesis, we can assume the result is true for the 
$(d-1)$-dimensional ring $\overline G$, which is the associated
graded ring for an $\fm$-primary ideal in a ring of dimension $d-1$. 
So we have  $\overline\omega=0$. In
other words $\omega \in (y_t^*)\omega_G$.
This means that  $\omega\in 
(y_1^*\cdots y_{t-1}^*)\omega_G \cap (y_t^*)\omega_G$.
Because the two element set 
$\{y_1^*\cdots y_{t-1}^*, \, y_t^*\}$ 
is an $\omega_G$-regular sequence, this means that 
$\omega \in (y_1^*\cdots y_t^*)\omega_G$ as wanted. 
This completes the proof  of the inductive step. 

\bigskip
It remains only to
prove the case where $d =1$. 
For this, we will invoke the careful description of $\omega_G$ in the
one dimensional case proved in Lemma \ref{dimonecase}. 
We will also need the  following two lemmas.
The first is a modification of
~\cite[Lemma 2.2]{CPU2}, which in turn is
inspired by \cite[Lemma 3.8]{HS}.
 The second is  a one-dimensional version of 
 \cite[Proposition 2.3]{Hy2} valid also in the non-Cohen-Macaulay case.

\begin{lem} 
\label{dimonelem}
Let $(A,\fm)$ be a local ring and let
$N \subset N' \subset M$ be  $A$-modules. Let $x,y\in A$, with 
 $x$ a nonzerodivisor on  $M$. Suppose that $yN\subset xN$. 
 Let $u_1,\ldots,u_m\in A$ be units of $A$ whose
images are distinct modulo $\fm$, and assume that $m \gg 0$
(specifically,    $m>\dim_k(N:_{N'}\fm)/N$, where $k = A/\fm$). Then 
\begin{equation*}x(N:_M\fm)\quad \bigcap  
\quad\bigcap_{i=1}^m(x+u_iy)(N:_{N'}\fm)
\subset xN'.\end{equation*}
\end{lem}

\begin{proof} Suppose that $w=xs=(x+u_1y)s_1= \ldots =(x+u_my)s_m$ where
$s\in N:_M\fm$ and
 $s_i\in (N:_{N'}\fm)$ for each  $i=1,\ldots,m$.
 For all $r\in N:_M\fm$,\, let 
$\overline r$ denote the class $r+N$ in $ N:_M\fm/N$. 
Also let $\overline a$ denote the class $a+\fm\in A/\fm$ for all $a\in A$.
The elements $\overline{u_1}\,\overline{s_1},\ldots,
\overline{u_m}\,\overline{s_m}\in (N:_{N'}\fm)/N$
 are  linearly dependent, since $m$ exceeds the dimension of this space.
By replacing $m$ by a possibly smaller positive integer, we will 
assume henceforth that every proper subset of the set  
$$\{\overline{u_1}\,\overline{s_1},\ldots,
\overline{u_m}\,\overline{s_m}\}$$
is linearly independent, but that the full set itself is dependent.
Clearly $m > 1$; otherwise, all $s_i$ are in $N$ and so 
$(x+u_iy)s_i \in xN \subset xN'$ as needed.

There exist units $\gl,\ldots,\gl_m\in A$
such that 
\begin{equation}\label{eq1}
\sum_{i=1}^m\overline{\gl_i}\,
\overline{u_i}\,\overline{s_i}=0.
\end{equation}
 This implies
\begin{equation*}\sum_{i=1}^m\gl_iu_is_i\in N,
\,\,\,{\text{ so that }} \,\,\,
(\sum_{i=1}^m\gl_iu_is_i)y\in yN\subset xN.\end{equation*}
Hence
\begin{equation*}(\sum_{i=1}^m\gl_iu_is_i)y=t x\end{equation*} for some $t\in N$.
Set \begin{equation*}\gl=\sum_{i=1}^m\gl_i.\end{equation*} Then
\begin{equation*}\gl s x=\sum_{i=1}^m\gl_is_i(x+u_iy)
=\sum_{i=1}^m\gl_is_ix+\sum_{i=1}^m\gl_is_iu_iy
=\sum_{i=1}^m\gl_is_ix+tx.\end{equation*} As $x$ is $M$-regular, this implies that
\begin{equation*}\gl s=\sum_{i=1}^m\gl_is_i+t.\end{equation*}

Now $\gl \in \fm$ would imply $\gl s\in N$ and so
\begin{equation*}\sum_{i=1}^m\overline{\gl_i}\,\overline{s_i}=0.
\end{equation*}
In this case, we can solve 
$$
\overline s_m = - \sum_{i= 1}^{m-1} 
\frac{\overline \lambda_i}{\overline \lambda_m} \overline s_i
$$
and plug into equation (\ref{eq1}) above.
Because the $\overline{u}_i$'s are pairwise distinct, this produces a non-trivial linear 
relation on the set 
$$\{\overline{u}_1\,\overline{s}_1,\ldots,
\overline{u}_m\,\overline{s}_{m-1}\},$$
a contradiction.
This forces 
 $\gl$ to be a unit in $A$,  and so  
\begin{equation*}s=\sum_{i=1}^m\gl^{-1}\gl_is_i+\gl^{-1}t\in N'.\end{equation*}
The proof of  Lemma \ref{dimonelem} is complete.   
\end{proof}

\begin{lem}\label{6.3dim1}
Let $I$ be an $\fm$-primary ideal in a 
local ring $(A, \fm)$  of dimension one.
Let $(x)$ be any minimal reduction of $I$. 
Then 
$$\Omega_0 = (x^t\omega_A:_{\omega_A} I^t)
$$
for all $t \gg 0$.
\end{lem}

\begin{proof}
Note that the blowup of $\spec A$ along $I$ is the 
affine scheme $Y = \Spec A[\frac{I}{x}]$. 
So the proper map $Y \la \Spec A$ corresponds to a 
finite  map of rings $A \rightarrow A[\frac{I}{x}]$.
Thus
$$
\Omega_0 = \omega_{\cO_Y} = \Hom_A(\cO_Y, \omega_A).
$$

To compute $\Omega_0$, 
note that without loss of generality, the $A$-module
generators for the ring $\cO_Y$ may be assumed of the
form $\frac{z}{x^t}$ where $z \in I^t$, for some fixed $t \gg 0$.
We claim that $\Omega_0 = (x^t\omega_A:_{\omega_A} I^t)$ for this
fixed $t$.

To check the inclusion  $\Omega_0 \subset (x^t\omega_A:_{\omega_A} I^t)$,
take any $f \in  \Hom_A(\cO_Y, \omega_A)$. 
Then the restriction of $f$ to $A$ is given by multiplication by 
$f(1) \in \omega_A$.
Because $f(\frac{z}{x^t}) \in \omega_A$ for any $z \in I^t$, 
 one readily verifies that $f(z) \in x^t\omega_A$ and so 
$f(1) I^t \in (x^t)\omega_A$. Then the 
map $f \mapsto f(1)$ gives the natural 
inclusion    $\Omega_0 \subset (x^t\omega_A:_{\omega_A} I^t).$

For the reverse inclusion, take any $w \in  (x^t\omega_A:_{\omega_A} I^t).$
For $z \in I^t$, we have $wz = x^tu$ for some $u \in \omega_A$.
Set $f_w(\frac{z}{x^t}) = u$. Using the fact that $x$ is a nonzerodivisor
on $\omega_A$,   one easily checks that the association
$w \mapsto f_w$ gives a well-defined
 injection   $(x^t\omega_A:_{\omega_A} I^t)
\subset \Hom_{A}(\cO_Y, \omega_A)$  inverse to the 
map in the previous paragraph. 
\end{proof}

We now prove the dimension one case of Theorem \ref{hard}, which will 
complete its proof, and hence the proof of the main technical theorem,
Theorem \ref{mainalg}. 

Let $d = 1$. 
Suppose that  $\omega\in [x^*\omega_G]_1$ for
all $x\in I$ which generate a minimal reduction of $I$
(of degree $n$ in the graded case).
As $x$ is not in any minimal prime of $A$, the element $x$ forms
 a system of
parameters for $A$ and so $x$ is  nonzerodivisor on the $A$-module
$\omega_A$.  
By Lemma~\ref{Glem} we know that $[\omega_G]_1
=\Omega_0/\Omega_1$. So we represent  $\omega$ by some
 for some $w\in
\Omega_0$, modulo $\Omega_1$.
We want to prove that $\omega=0$, which is the same thing as
proving that $w \in \Omega_1$. 
 
Because the module $[\omega_G]_1 = \Omega_0/\Omega_1$ is an Artinian 
$A$-module,
 we might as well  assume that
$\fm\omega=0$. In other words, we can assume
that $w\in \Omega_1:_{\Omega_0}\fm$.
Again by Lemma~\ref{Glem}, we know that  $[\omega_G]_0
=(\Omega_{-1}\cap \omega_A)/\Omega_0$. So, using also that 
$\Omega_1 = x\Omega_0$ by Lemma \ref{lem1} we can write 
 $w=sx$ for some $s\in \Omega_{-1}\cap \omega_A$. 
Since
 $$
\fm s x=\fm w \in \Omega_1 = x\Omega_0,
$$
we have that 
  $$
\fm s\in \Omega_0,
$$
as $x$ is a nonzerodivisor on $\omega_A$.
 This means that $s \in \Omega_0:_{\omega_A}\fm$. 
Our goal is to show $s \in \Omega_0$. This will complete the proof
since then $w = x s \in I\Omega_0 \subset \Omega_1$, whence the class of 
$w$ in $\omega_G$ is zero.

\smallskip
To achieve our goal of showing that 
 $s \in \Omega_0$,  we invoke Lemma \ref{6.3dim1},
which guarantees that it is enough to show that 
$s \in x^t\omega_A:I^t$ for $t \gg 0$.
 For this, note that 
it is enough to show that 
\begin{claim}\label{claim1}
For any $y \in I$ (homogeneous in the graded case),  
  $s\in 
(x^t\omega_A :_{\omega_A} y^t)$ for all $t\gg 0$.  
\end{claim}
Indeed,  because $\bbQ \subset A$, the ideal generated by the
$t$-th powers of the 
 elements of an ideal $I$ (homogeneous elements of degree $n$ when $I$ is
generated in degree $n$)
  is simply $I^t$.
This follows easily from the identity
\begin{equation}\label{identity1}
t!X_1\cdots X_t
=\sum_{1\le k\le t}\sum_{1\le i_1<\ldots<i_k\le t}(-1)^{t-k}
(X_{i_1}+\ldots+X_{i_k})^t.
\end{equation}

\medskip
To prove Claim \ref{claim1},
 take any  $y\in I$ (of degree $n$ in the graded case).
We wish to apply Lemma \ref{dimonelem} to the $A$-modules
$$
\Omega_0 \subset N_t' \subset \omega_A,
$$
where $N_t' = (x^t\omega_A :_{\omega_A} y^t)$.
Note here  that $\Omega_0 \subset N_t'$ for all $t$, since
$$
y^t\Omega_0 \subset I^t\Omega_0 \subset \Omega_t = x^t\Omega_0 
\subset x^t\omega_A.
$$

To this end, choose distinct units $u_1, \dots,  u_m$ (of degree zero in the
graded case), 
with $m \gg 0$, and so that each $x + u_iy$ generates a reduction for $I$.
Then
\begin{equation}\label{identity}
w=xs=(x+u_1y)s_1=\ldots =(x+u_my)s_m
\end{equation} for some $s_1,\ldots,s_m\in \Omega_0:_{\omega_A}\fm$.  
Therefore, assuming for the moment  that the elements $s_i$ are also 
in $N_t'$,
then  $xs$ is in 
$$
x(\Omega_0:_{\omega_A} \fm) \cap [(x+u_1y)(\Omega_0:_{N_t'}\fm)]
 \cap \dots \cap 
[(x+ u_my)(\Omega_0:_{N_t'}\fm)],$$
and so applying Lemma \ref{dimonelem}, we have
$$
xs \in x N_t'.$$
Because $x$ is a nonzerodivisor on $\omega_A$ (and so on any submodule),
we see that $s \in N_t' = (x^t\omega_A :_{\omega_A} y^t)$, which is
precisely what we needed to show.

It remains to show that each $s_i$ is in $ N_t' = 
(x^t\omega_A :_{\omega_A} y^t)$
for $t \gg 0$. We in fact will show this for all $t \geq 0$, using 
induction on 
 $t$.
 If $t=0$, the $s_i$ are trivially in $N'_0$.
Assume then that 
\begin{equation*}
s_1,\ldots,s_m\in N_{t-1}' =  (x^{t-1}\omega_A :_{\omega_A} y^{t-1}).
\end{equation*}
So by the argument above using Lemma \ref{dimonelem},
we have also that $s \in N'_{t-1}$. 
Since 
$$xs = (x + u_iy)s_i,$$
we have
$$ 
ys_i = u_i^{-1}x(s - s_i) \in x N_{t-1}' = 
x(x^{t-1}\omega_A :_{\omega_A} y^{t-1}) 
\subset x^{t}\omega_A :_{\omega_A} y^{t-1}.
$$
But then for each $i = 1, \dots, m$, we have 
$$s_i \in  (x^{t}\omega_A :_{\omega_A} y^{t}) = N'_t,
$$
as we sought to show. 

\medskip
 This completes the proof of Claim \ref{claim1}, and hence the proof of 
Theorem \ref{hard} in the dimension one case, and the proof of our main 
technical theorem, Theorem
\ref{mainalg}.
\end{proof}

\begin{rem}\label{d=2} In  Theorem \ref{mainalg},
as well as in its later applications in Section 5, 
the assumption that  $A$ contains the set of rational numbers
is unnecessary  
in case the
ideal $I$ has reduction number at most one. Note that the
reduction number does not increase as we reduce to the one-dimensional
case. Now it is easy to see that in Lemma \ref{6.3dim1}, one can take 
$t$ to be the reduction number, that is, 
$\Omega_0 = (x^t\omega_A:_{\omega_A} I^t) = (x^r\omega_A:_{\omega_A} I^r)$
where $r$ is the reduction number of $I$. So when the reduction number of
$I$ is at most one, it suffices to prove  Claim \ref{claim1} for the 
case $t = 1$. This follows in the same way without making use of Identity
(\ref{identity1}).
\end{rem}

\medskip

\section{Brian\c con-Skoda Type results.}

In this section, we identify general conditions under which 
the hypothesis of our Main Technical Theorem \ref{mainalg}
are satisfied. The  results of this section will allow us to 
apply the main technical theorem proved in the preceding section to
deduce our main results in the following section.

\subsection{The Brian\c con-Skoda Theorem and related properties of adjoints.}
\label{BStype}
The next result is 
analogous to the ``Brian\c con-Skoda Theorem with adjoints'' proved by Lipman
for regular schemes in \cite{L3}. See also the ``Skoda Theorem''
discussed  in \cite{laz}.

\begin{lem}
\label{BS}
Let $(A,\fm)$ be a local ring  of dimension $d \geq 1$,
and let $I\subset A$ be an   ideal of positive height.
Assume that $\omega_Y$ is $m$-regular in the sense of 
Castelnuovo-Mumford, meaning that
$$
H^i(Y, \omega_Y(m-i))= 0
$$ for all $ i > 0$. Then  
 $\Omega_{n+1}=J\Omega_n$ for
all $n\ge m$, where $J$ is any reduction of $I$.
\end{lem}

\begin{proof}
 This follows from a standard  argument.
 As $I^{r+1}=JI^r$, we have $J\cO_Y=I\cO_Y=\cO_Y(1)$. Generators of $J$ 
therefore give rise to global sections which generate the sheaf $\cO_Y(1)$.
These give rise to an exact Koszul complex. One can then 
argue as in~\cite[p.~747]{L3} (or as below in the proof of Proposition
\ref{intersectionlem1})
 to conclude that $\Omega_{n+1}=J\Omega_n$ 
for all $n \geq m$ under the stated vanishing conditions.
\end{proof}

\begin{cor}\label{BS2}
If $I$ is $\fm$-primary, and  $A$ and $R = A[It]$ are Cohen-Macaulay, then  
 $$\Omega_{n}=J\Omega_{n-1}$$ 
for all $n \geq d = \dim A.$ 
\end{cor}

\begin{proof}
By Proposition \ref{BS}, it suffices to check that 
 $\omega_Y$ is $(d-1)$ regular.
But the
  the vanishing of 
 $H^i(Y, \omega_Y(d-1-i) )$  for positive $i$ is
the same as the vanishing of $H^{i}_Z(Y, \cO_Y(-i+1))$ for 
$i \leq d-1$. This vanishing
follows easily from the Sancho de Salas sequence for $R$, as
 recorded in
\ref{SSfacts}.
\end{proof}

\begin{rem}
Corollary \ref{BS} holds also for ideals $I$ that have a reduction 
generated by regular sequence, with $d$ now the height of $I$. 
Furthermore, 
 because $I\cO_Y = \cO_Y(1)$ is ample for the map $\pi$, 
we have  $H^i(Y,\omega_Y(m-i))=0$ 
for all $i>0$ and  for sufficiently large $m$.  So for large enough $n$,
we always have $\Omega_{n+1} = J\Omega_n$, where $J$ is any reduction of $I$.
\end{rem}

\bigskip
The next result
ensures that  Theorem \ref{mainalg} can be applied to 
 some interesting cases.
Note that when $n \geq d$,
this statement
collapses to the Brian\c con-Skoda Theorem above in \ref{BS2}; the following
result thus informs us also of what goes on for smaller indices $n$.

\begin{prop} 
\label{intersectionlem1} Let  $I$ be an $\fm$-primary ideal in a
Cohen-Macaulay  ring $(A, \fm)$ of dimension  
$d> 0$, and let $J$ be any minimal reduction of $I$.
Let $R$ denote the Rees ring $A[It]$ of $A$ with respect to $I$, and assume 
that $R$ is Cohen-Macaulay. 
Then 
$$
J\omega_A\cap \Omega_n=J\Omega_{n-1}
$$ for all
 $n\in \bbZ$.
\end{prop}

\begin{proof} 
The case $d =1$ is degenerate. 
 The condition that $R$ is Cohen-Macaulay forces
$I$ to be principle (see e. ~g. \cite[Cor 25.2]{HIO}),
 so $I = J = (x)$.
Then  $Y = \spec A$ and $\Omega_n = x^n\omega_A$ for all $n$.
This makes the statement obvious. 

Assume $d > 1$.
Because $A$ and $R$ are both Cohen-Macaulay,  
$\Omega_n = \omega_A$ for non-positive $n$; see Lemma 
\ref{stab} and \ref{SSfacts}.
So 
the statement is trivial for negative $n$.  Also, when $n \geq d$, 
the statement follows from Corollary \ref{BS2}. 
It remains to consider the case $1 \leq n < d$.

\medskip 
Let us first 
consider the weaker statement that
\begin{equation*}
J\Omega_{n-2}
\cap \Omega_n=J\Omega_{n-1}
\end{equation*}
for all $n<d$.  
To prove this, it is sufficient to prove that the natural map
\begin{equation*}
\Omega_n/J\Omega_{n-1}\lrarrow \Omega_{n-1}/J\Omega_{n-2}
\end{equation*} 
is injective.

Because $I^{r+1}=JI^r$ for some $r\ge 0$, we have
 $I{\cO}_Y=J{\cO}_Y$. Fix generators 
 $\{x_1,\ldots, x_d\}$ for $J$ as an ideal of  $A$.
These elements 
give
rise to generating global sections of $I{\cO}_Y$,
and hence a surjection
$\bigoplus_{i=1}^d{\cO}_Y \lrarrow I{\cO}_Y$,
 and so also a 
surjection
$$\cF := \bigoplus_{i=1}^d {\cO}_Y(-1)
\overset{\begin{pmatrix}x_1 \\ \vdots \\ x_d\end{pmatrix}}{\la} 
\cO_Y.
 $$
This means that the Koszul complex 
\begin{equation*}
0\lrarrow \wedge^d{\cF}\otimes
\omega_Y(n) \lrarrow \dots \lrarrow\wedge^1 {\cF}\otimes \omega_Y(n)\lrarrow
\omega_Y(n)\lrarrow 0 \end{equation*}
 is exact, where 
\begin{equation*}
\wedge^j {\cF}\otimes \omega_Y(n)=(\omega_Y(n-j))^{\oplus {d
\choose j}},
 \end{equation*}  
for all $j = 0, 1, \dots, d.$
We now split this complex into $d-1$ short exact sequences
\begin{equation*}0\lrarrow {\cK}_j \lrarrow \wedge^j {\cF}\otimes \omega_Y(n)
 \lrarrow 
{\cK}_{j-1}\lrarrow 0,
\end{equation*} 
where 
${\cK}_0=\omega_Y(n)$ and ${\cK}_{d-1}=\omega_Y(n-d)$.
The corresponding long exact sequences of cohomology give the exact sequences
\begin{equation*}
H^{j-1}(Y,\wedge^j {\cF}\otimes \omega_Y(n))\lrarrow H^{j-1}
(Y,{\cK}_{j-1})\lrarrow H^j(Y,{\cK}_j)
\end{equation*} 
for each $j = 1, 2, \dots, d-1$.
When $j=1$ we get the sequence 
\begin{equation*} \bigoplus_{i=1}^d \Omega_{n-1}
\overset{\begin{pmatrix}x_1 \\ \vdots \\ x_d\end{pmatrix}}{\la}
 \Omega_n
\lrarrow H^1(Y,\cK_1) 
\end{equation*}
 which gives an injection
\begin{equation*}  \Omega_n/J\Omega_{n-1}
\hookrightarrow  H^1(Y,\cK_1).\end{equation*} 
Note that 
 $H^{j-1}(Y,\omega_Y(n-j))=0$ for $1<j\le d-1$;
this 
follows easily from the 
 Sancho de Salas sequence (see  \ref{SSfacts}),
taking into consideration the abundant vanishing afforded 
because $A$ and $R$ are Cohen-Macaulay. 
Thus  
we obtain injections
\begin{equation*} H^{j-1}(Y,{\cK}_{j-1})\hookrightarrow H^j(Y,{\cK}_j)
\end{equation*} 
and finally, an injection
\begin{equation*}  \Omega_n/J\Omega_{n-1}
\hookrightarrow H^{d-1}(Y,\omega_Y(n-d)).\end{equation*}

The inclusion $I^n\omega_Y\subset I^{n-1}\omega_Y$  induces 
a homomorphism of complexes 
\begin{equation*}
\begin{array}{ccccccccc}
0&\lrarrow& \wedge^d{\cF}\otimes\omega_Y(n)& \lrarrow 
\dots \lrarrow&\wedge^1 {\cF}\otimes \omega_Y(n)&\lrarrow&
\omega_Y(n)\\
\Big\downarrow& &\Big\downarrow&  &\Big\downarrow& &\Big\downarrow\\
0&\lrarrow& \wedge^d{\cF}\otimes\omega_Y(n-1) &\lrarrow \dots \lrarrow&\wedge^1 {\cF}\otimes \omega_Y(n-1)&\lrarrow&
\omega_Y(n-1),
\end{array}
\end{equation*}
which in turn induces a 
  commutative diagram 
\begin{equation*}
\begin{CD}
\Omega_n/J\Omega_{n-1}
@>>>H^{d-1}(Y,\omega_Y(n-d))\\
@VVV @VVV\\
\Omega_{n-1}/J\Omega_{n-2}
@>>>H^{d-1}(Y,\omega_Y(n-1-d))\\
\end{CD}
\end{equation*}
Thus we can get the desired injectivity of
$\Omega_n/J\Omega_{n-1}\lrarrow \Omega_{n-1}/J\Omega_{n-2}$
by proving 
injectivity for   the homomorphism
\begin{equation*}
H^{d-1}(Y,\omega_Y(n-d))\lrarrow H^{d-1}(Y,\omega_Y(n-1-d)).\end{equation*}

To this end, consider
the adjunction sequence for $E \subset Y$, 
\begin{equation*}0\lrarrow \omega_Y \lrarrow \omega_Y(-1)
\lrarrow \omega_E\lrarrow 0. 
\end{equation*}
Tensor with the flat module $I^{n-d}\cO_Y$ to get a short exact 
sequence
\begin{equation*}
0\lrarrow \omega_Y(n-d) \lrarrow \omega_Y(n-1-d) \lrarrow \omega_E(n-d) \lrarrow 0.
\end{equation*} 
Taking cohomology, we get the exact sequence
\begin{equation*}
H^{d-2}(E,\omega_E(n-d)) \lrarrow H^{d-1}(Y,\omega_Y(n-d))
 \lrarrow H^{d-1}(Y,\omega_Y(n-1-d)).\end{equation*}
Thus it is sufficient   to show that $H^{d-2}(E,\omega_E(n-d))=0$. 
By  duality applied to the map $E \la \spec A/I$,
 it is enough to show
$ H^1(E,\cO_E(d-n)) = 0.$ 
But this follows easily from the
 Sancho de Salas
sequence for the graded ring $G$:
$$
H_{\fm}^1(G_{d-n})\lrarrow H^1(E,\cO_E(d-n)) \lrarrow [H_{\fM_G}^2(G)]_{d-n}.
$$
Indeed, because $n < d$, the module $ H_{\fm}^1(G_{d-n}) = 
H^1_m(I^{d-n}/I^{d-n+1})$ vanishes (as the module $I^{d-n}/I^{d-n+1}$
has zero-dimensional support).
Now 
when 
 both   $R$ and $A$ 
are Cohen-Macaulay, the associated graded ring $G$ is
Cohen-Macaulay with negative $a$-invariant (see \ref{SSfacts}), so
the module $[H_{\fM_G}^2(G)]_{d-n}$ certainly vanishes,
as $d - n > 0$.

\medskip
We have now shown that $$
J\Omega_{n-2}\cap \Omega_n = J \Omega_{n-1}
$$ for all $n <d$. To complete the proof, note
$J\Omega_{n-1} \subset J\omega_A \cap \Omega_n$ for all $n \geq 1.$
So we need to show the reverse inclusion. We will do this by induction
on $n$, starting from $n = 0$. 

For $n = 0$, we verify that $J\omega_A \cap \Omega_0 \subset J\Omega_{-1}$.
Indeed, because the $a$-invariant of $G$ is non-positive, we know that
$\omega_{A} \subset \Omega_{-1}$. So  $J\omega_{A} \subset J\Omega_{-1}$,
and of course then $J\omega_{A} \cap \Omega_0 \subset J\Omega_{-1}$.

Now assume that $n >0$ and the inclusion has been proved  for smaller
indices. Take
 $x\in J\omega_A \cap \Omega_n$. Then certainly $x \in  
J\omega_A \cap \Omega_{n-1}$, which is $J\Omega_{n-2}$ by the induction
hypothesis.
So  $x\in J\Omega_{n-2}\cap \Omega_n$. But by
the weaker statement proved above, this implies that $x\in J\Omega_{n-1}$. 
The proof is complete.
\end{proof}

\subsection{The case where $A$ is not necessarily Cohen-Macaulay.}

The following proposition offers an even more general setting in 
which the hypothesis of our main technical theorem are satisfied.
Its proof is decidedly less elementary than the
argument we have already made for Proposition \ref{intersectionlem1}
(which is why we have included a separate proof for  \ref{intersectionlem1}).

\begin{prop} 
\label{intersectionlem2} 
Let $(A,\fm)$ be a local ring of dimension $d> 1$. Let $I\subset A$ 
be an $\fm$-primary ideal such that the irrelevant ideal  
of the Rees ring $A[It]$ is a Cohen-Macaulay $A[It]$-module. 
 Then, for any reduction $J$ of $I$, 
 $$J\omega_A\cap \Omega_n=J\Omega_{n-1}$$ for all $n\in \bbZ$.
\end{prop}

\begin{proof} Let $R^+$ denote the irrelevant ideal of $R = A[It]$. 
Because $R^+$ is Cohen-Macaulay, the $a$-invariant of the associated graded 
ring 
 $G$ is non-positive,
as one checks by  looking at the Sancho de Salas sequence. 
So from Lemma \ref{stab}, we have that $\Omega_n = \omega_A$ for $n < 0$.
Thus statement is trivial for $n \leq -1$. Also the case $n = 0$ is clear,
since $J\omega_A \cap \Omega_0 = J\Omega_{-1} \cap \Omega_0 = J\Omega_{-1}.$
For $n >0$ we proceed by induction on $n$ as in the proof of Proposition
\ref{intersectionlem1}. As in that argument, 
it is sufficient to prove
 the weaker statement that
\begin{equation*}J\Omega_{n-2}\cap \Omega_n=J\Omega_{n-1}\end{equation*}for all $n\ge 1$.  

Choose $N \gg 0$ such that $\Omega_{N+1}= {J} \Omega_N$.
Write 
\begin{equation*}
l_A(\Omega_{-1}/ \Omega_{N+1})= l_A(\Omega_{-1}/ \Omega_0)+\sum_{n=0}^{N}l_A(\Omega_{n}/\Omega_{n+1}),
\end{equation*} 
where $l_A(M)$ denotes the length of an $A$-module $M$, 
and 
\begin{equation*}
l_A(\Omega_{-1}/J\Omega_N)=l_A(\Omega_{-1}/J\Omega_{-1}) + 
\sum_{n=0}^Nl_A(J\Omega_{n-1}/J\Omega_n).\end{equation*}
Then 
\begin{equation*}
l_A(\Omega_{-1}/{J}\Omega_{-1})= 
l_A(\Omega_{-1}/\Omega_0)+ \sum_{n=0}^{\infty}(l_A(\Omega_{n}/\Omega_{n 
+1})- l_A({J} \Omega_{n-1}/{J} \Omega_{n})).\end{equation*}
Consider the $G$-module $W=\bigoplus_{n\ge 0}\Omega_{n-1}/\Omega_n$. 
Fix generators  $\{x_1,\ldots, x_d\}$ for $J$ as an ideal of  $A$. Let $J^*$ denote the ideal
$(x_1^*,\ldots,x_d^*)\subset G$. 
Because 
\begin{align*}
\begin{split}
l_G(W/J^*)&=l_A(\Omega_{-1}/\Omega_0)+
 \sum_{n=1}^{ \infty} l_A(\Omega_{n-1}/J\Omega_{n -2} + 
\Omega_n)\\
&=l_A(\Omega_{-1}/\Omega_0)+ \sum_{n=1}^{ 
\infty} (l_A(\Omega_{n-1}/\Omega_{n})-l_A({J}\Omega_{n -2}+ \Omega_n/\Omega_n))
\\
&= l_A(\Omega_{-1}/\Omega_0)+ \sum_{n=1}^{ \infty} 
(l_A(\Omega_{n-1}/\Omega_n)-l_A(J\Omega_{n -2}/J\Omega_{n-2}\cap 
\Omega_n))\\
&= l_A(\Omega_{-1}/\Omega_0)+ \sum_{n=0}^{ 
\infty} (l_A(\Omega_{n}/\Omega_{n+1})-l_A({J}\Omega_{n-1}/{J}\Omega_{n-1}\cap 
\Omega_{n+1})),
\end{split}
\end{align*} 
we now obtain 
\begin{align*}
\begin{split}
l_G(W / J^*)-l_A(\Omega_{-1}/J\Omega_{-1})&= 
\sum_{n=0}^{\infty} (
l_A({J} \Omega_{n-1}/{J} \Omega_n))
-l_A({J}\Omega_{n-1}/{J}\Omega_{n-1}\cap \Omega_{n+1}))\\
&=
\sum_{n=0}^{\infty} l_A(J\Omega_{n-1}\cap \Omega_{n+1}/ J\Omega_n).
\end{split}
\end{align*}
In order to prove our claim we thus have to 
show that $l_G(W / J^*)-l_A(\Omega_{-1}/J\Omega_{-1})=0$.
For this,  recall the notion of the $\bbI$-invariant of a 
graded module (see e.~g.~\cite[p.~6]{T}). 
Let $B$ be a graded ring defined over
 a local ring, and let $N$ be a graded
$B$-module with $r=\dim N$. Then the 
$\bbI$-invariant  
\begin{equation*}\bbI(N)=\sum_{i=0}^{r-1}{r-1\choose i}l_B(H_{\fN}^i(N))\end{equation*} 
where $\fN$ denotes the homogeneous maximal ideal of $B$.
It is 
a general fact that 
if $(y_1,\ldots,y_r)$ is any homogeneous system of parameters for $B/\ann{N}$,
 then always
$l_A(N/(y_1,\ldots,y_r)N)-e(y_1,\ldots,y_r;N)\le \bbI(N)$, where
$e$ denotes the multiplicity of $N$ with respect to $(y_1, \dots, y_r)$.
 If the equality holds, then the system
of parameters $(y_1,\ldots,y_r)$ is called standard. 

By definition the multiplicity $e(J^*;W)=e(G^+;W)$ is equal to d! times
the leading coefficient of the numerical polynomial $l_A(W_n)$ where $n\gg 0$. But for $n\gg 0$,
\begin{equation*}W_n=\Omega_{n-1}/\Omega_n=J^{n-1-N}\Omega_N/J^{n-N}\Omega_N\end{equation*}
showing that $e(J^*;W)=e(J,\Omega_N)$. 
For any $\fp\in \Spec A$, $\fp\not\supset I$, we clearly 
have $(\Omega_n)_{\fp}=\omega_{A_p}$ for all $n\in \bbZ$.
Therefore  \begin{align*}
\begin{split}
e(J,\Omega_N)&=
\sum_{\fp\in \Min A, \dim A/\fp=d}l_{A_{\fp}}
((\Omega_N)_{\fp})e(J+\fp/\fp;A/\fp)\\
&=
\sum_{\fp\in \Min A, \dim A/\fp=d}l_{A_{\fp}}
(\omega_{A_{\fp}})e(J+\fp/\fp;A/\fp)\\
&=
\sum_{\fp\in \Min A, \dim A/\fp=d}l_{A_{\fp}}
((\Omega_{-1})_{\fp})e(J+\fp/\fp;A/\fp)\\
&=
e(J;\Omega_{-1}).
\end{split}
\end{align*}
So $e({J}^*;W)=e(J;\Omega_{-1})$. 

Set $\Omega=\bigoplus_{n\ge -1}\Omega_n$
and $\Omega'=\bigoplus_{n\ge 0}\Omega_n$. Observe that 
$\Omega'=\omega_{R^+}$ where
$\omega_{R^+}$ denotes the canonical module of the $R$-module $R^+$,
that is, $\Omega'$ is the graded dual of  the top local cohomology module of
$R^+$ with supports in the unique homogeneous maximal ideal.  
Indeed, a look at the Sancho de Salas sequence
\begin{equation*}
0 \lrarrow  H_{\fm}^d(R^+_n) \lrarrow H_E^d(X,\cO_X(n))\lrarrow
 [H^{d+1}(R^+)]_n\lrarrow 0
\end{equation*} 
shows that there is an isomorphism
$H_E^d(X,\cO_X(n))=[H^{d+1}(R^+)]_n$, for all $n\le 0$.
 Dualizing, we then see that
$\Omega_n=[\omega_{R^+}]_n$ for all $n\ge 0$. On the other hand, by 
considering the long exact sequence of 
cohomology
corresponding
 to the exact sequence \begin{equation*}0 \lrarrow R^+\lrarrow R\lrarrow A
\lrarrow 0\end{equation*}
and taking into account that $a(R)=-1$,
we get $[H^{d+1}(R^+)]_n=[H^{d+1}(R)]_n=0$
 for $n>0$ so that $[\omega_{R^+}]_n=0$ when $n<0$.
By ~\cite[Satz 3.2.2]{S}, it follows that 
 $\Omega'$ is Cohen-Macaulay.
 By means of the long exact sequences of cohomology corresponding 
to the exact sequences \begin{equation*}
0\lrarrow \Omega'\lrarrow \Omega  \lrarrow \Omega_{-1}(1)\lrarrow 0
\quad\text{and}\quad 0\lrarrow \Omega'\lrarrow 
\Omega(-1) \lrarrow W\lrarrow 0 \end{equation*} 
one 
then
easily checks that $H_{\fM}^i(W)=H_{\fm}^i(\Omega_{-1})$ for
all $0\le i<d$. Therefore $\bbI(W)=\bbI(\Omega_{-1})$. As $a(G)\le 0$, 
Lemma \ref{stab} implies that $\Omega_{-1}=\omega_A$. 
Using~\cite[Theorem 1.1, Appendix]{GY} and~\cite[Corollary 6.18]{GY}
 we know that $(x_1,\ldots,x_d)$ is a 
standard system of parameters for $A$. Moreover, by~\cite[Theorem 3.17]{GY},
then  
$(x_1,\ldots,x_d)$ is
standard also for $\omega_A$. 
This implies that
$\bbI(\omega_A)=l_A(\omega_A/J\omega_A)-e(J;\omega_A)$. Therefore
\begin{align*}
\begin{split}
l_G(W / J^*)-l_A(\Omega_{-1}/J\Omega_{-1})&=
(l_G(W / J^*)-e(J^*;W))\\
&-(l_A(\Omega_{-1}/J\Omega_{-1})-e(J;\Omega_{-1}))\\
&\le \bbI(W)-\bbI(\omega_A)=0
\end{split}
\end{align*}
as wanted. The proof is complete.
\end{proof}

\section{The main local algebraic results}
In this section, we prove our main results
about the core in a local ring, including the theorems
relating core, adjoint, and coefficients ideals,
the ``local'' version of Kawamata's Conjecture, and a
formula for core conjectured in \cite{CPU2}.  All are deduced from the main
technical theorem, Theorem \ref{mainalg}, using the Brian\c con-Skoda results
of Section 4. Further corollaries for graded rings
appear at the end of Section 6.
 
\smallskip

\subsection{Formulas for Core in the Cohen-Macaulay case.}
Recall that an {\it equimultiple} ideal in a local ring 
$(A, \fm)$ is an ideal whose height equals its analytic spread 
(see \cite[p58]{HIO}). When the ring $A$ has an infinite residue field,
an equimultiple ideal is precisely an ideal having a reduction generated
by part of a system of parameters. Every $\fm$-primary ideal in a local
ring is equimultiple. 

\begin{cor}\label{hswan}
Let $(A, \fm)$ be a Cohen-Macaulay local ring containing a
field of characteristic zero, and let $I$ be an equimultiple ideal of height
 $h$ whose  Rees ring $A[It]$ is  Cohen-Macaulay. 
Then 
$$
 \core (I) =  H^0(Y, I^h\omega_Y):_A \omega_A = J^{r+1}:_A I^r
$$
where, as always, $ H^0(Y, I^h\omega_Y)$ is considered as a 
submodule of $\omega_A$ via the trace map, and $J$ is any reduction of $I$.
\end{cor}

\begin{rem}
In Corollary \ref{hswan}, if the height $h$ is two (or less), 
then the assumption that $A$ contains the rational numbers is not
needed. See Remark \ref{d=2}.
 The case where $h$ is height one is trivial, because our
assumption on the Rees ring forces $I$ to be principal. 
\end{rem}

\begin{rem} The formula $\core(I)=J^{r+1}:I^r$  is conjectured 
in~\cite[Conjecture 5.1]{CPU2} under more general  hypothesis.
\end{rem}

Corollary \ref{hswan} 
 follows easily from the following theorem, which generalizes 
our main technical theorem to ideals that may not be $\fm$-primary. 

\begin{thm} 
\label{corethm}
Let $(A,\fm)$ be a Cohen-Macaulay local ring containing the set of rational 
numbers.
Let $I\subset A$ be an equimultiple ideal of positive height $h$ such that 
the corresponding Rees ring $A[It]$ is Cohen-Macaulay.  
Then  
$$\core(I\omega_A)= \Gamma(Y, I^h\omega_Y) 
$$ as submodules  of $\omega_A$,
where $Y = \proj A[It]$.
\end{thm}

\begin{rem}
In fact,  $\Omega_h \subset \core{(I\omega_A)} $ for equimultiple
ideals $I$ in a Cohen-Macaulay local ring $A$ without any assumption on the
Rees ring. The point is to prove the reverse inclusion. 
As we'll see in the proof, the reverse inclusion holds even when 
it is assumed only that the irrelevant ideal of $R$ is Cohen-Macaulay
(even when  $A$ is not). 
\end{rem}

To see that that Theorem \ref{corethm} implies Corollary \ref{hswan},
we need the following lemma.

\begin{lem}\label{lem6}
Let $(A, \fm)$ be a Cohen-Macaulay local ring, and let $I$ be an 
equimultiple ideal of height $h$.
Let $J$ be any minimal reduction of $I$. Then 
$$
\Omega_h :_A \omega_A = J^{r+1}:_A I^r
$$
for any integer $r$ such that $I^{r+1} = JI^r$.
In particular,  $J^{r+1}:I^r$ is independent of the choice
of reduction $J$.
\end{lem}

\begin{proof}
We first note that 
$$
\Omega_n=J^{n-h+r+1}\omega_A:_{\omega_A}I^r
$$ 
for all $n\ge 0$.
Indeed,  set $X=\Proj A[Jt]$ and $Y=\Proj A[It]$.
 Since $J$ is generated
by a regular sequence,  $\omega_{A[Jt]}$ has the 
expected form:
 $$[\omega_{A[Jt]}]_n=\Gamma(X,J^n\omega_X)=J^{n-h+1}\omega_A$$ 
for all $n\ge 1$ (see e.g.~\cite[p.~142]{Va}). 
The canonical sheaf for $Y$ can therefore be computed from the
canonical sheaf for $X$ via the finite map $Y \la X$ induced by the
inclusion of the Rees rings $A[Jt] \hookrightarrow A[It]$.
Working this out, we arrive at 
\begin{equation}\label{eq4} 
\Gamma(Y,I^n\omega_Y)=\Gamma(X,J^{n+r}\omega_X):_{\omega_A} I^r
\end{equation} as submodules of $\omega_A$ for all $n\ge 0$. 
For details, see 
\cite[Proposition 2.3]{Hy2}. 

Now we note that because 
$J$ is generated by a regular sequence,
\begin{equation}\label{eq5}
J^n\omega_A:_A\omega_A=J^n
\end{equation}
 for all $n\ge 0$. 
This can be proved by induction on $n$.
When   $n=1$, this follows 
from the fact that $\omega_{A/J}=\omega_A/
J\omega_A$ is a faithful $A/J$-module. 
Suppose then that $n>1$. Let $x\in J^n\omega_A:\omega_A$. By the induction
hypothesis, we know that $x\in J^{n-1}$.  There thus exists a form 
$F\in A[t_1,\ldots,t_r]$ of degree $n-1$ such that 
$x=F(a_1,\ldots, a_r)$. As $F(a_1,\dots,a_r)\omega_A\in J^n\omega_A$
 and $(a_1,\ldots,a_r)$ is an $\omega_A$-regular sequence,
every coefficient of $F$ must lie in $J\omega_A:\omega_A=J$.
  Hence $x\in J^n$.

\smallskip
Finally, to see that 
$$
\Omega_h:_A \omega_A = J^{r+1}:_A I^r,
$$
we simply compute
\begin{equation*}
\Omega_h:_A\omega_A
=(J^{r+1}\omega_A:_{\omega_A}I^r):_A\omega_A
=(J^{r+1}\omega_A:_A\omega_A):_AI^r
=J^{r+1}:I^r.
\end{equation*}
Here, the first equality follows from (\ref{eq4}) above,
and the last equality follows from (\ref{eq5}).
\end{proof}
\smallskip

\begin{proof}[Proof  that Theorem \ref{corethm} implies Corollary  \ref{hswan}]
Say that $y \in \core(I)$. Then
 $$y\omega_A \subset \core( I\omega_A) = \Omega_h$$ by
 Theorem \ref{corethm}, so $y \in \Omega_h:_A \omega_A = J^{r+1}:_A I^r$. 
Conversely, say that 
$$
y \in \Omega_h:_A\omega_A = J^{r+1}:_AI^r.
$$
Then $y \in  J^{r+1}:_AJ^r = J$ for every reduction $J$ of $I$, so 
$y \in \core(I)$. The proof of  \ref{hswan}
is  complete. 
\end{proof}

\smallskip

\begin{proof}[Proof of Theorem \ref{corethm}]
Assume  that $I$ is equimultiple.  As usual, write $\Omega_h$
for $\Gamma(Y, I^h\omega_Y)$, considered as a submodule of $\omega_A$.
To see that $\Omega_h \subset \core(I\omega_A),$
recall that $\Omega_h = J^{r+1}\omega_A :_{\omega_A} I^r$
 (as in  the proof of Lemma \ref{lem6}).
Thus $\Omega_h \subset  J^{r+1}\omega_A :_{\omega_A} J^r \subset J\omega_A$,
since $J$ is generated by a regular sequence.

It remains to show that $\core(I\omega_A)\subset \Omega_h$. 
Fix
elements $x_1,\ldots,x_{d-h}\in A$ such that
$(x_1+I,\ldots,x_{d-h}+I)$ is a system of parameters on $A/I$ and a 
regular sequence 
of degree
 zero elements on $G$, the 
associated graded ring of $A$ with respect to $I$
 (see e.g.~\cite[Proposition 10.24 and its proof]{HIO}). 
Then also $(x_1^t+I,\ldots,x_{d-h}^t+I)$ is a regular sequence on $G$. 
Set $\overline A=A/(x_1^t,\ldots,x_{d-h}^t)$ and let 
$\overline I$ denote the image of $I$ in $A$.
 Then $\overline I$ is 
an $\overline\fm$-primary ideal of $\overline A$.
 Moreover, the corresponding associated graded ring
$$
\overline G = \gr_{\overline I}{\overline A} :=
\frac{\overline A}{\overline I} \oplus 
\frac{\overline I}{{\overline I}^2} \oplus \dots
$$ is easily seen to be isomorphic to 
$$
G/(x_1^t+I,\ldots,x_{d-h}^t+I),$$ 
and hence 
Cohen-Macaulay.
(For example, one can use the fact that 
 $(x_1^t,\ldots,x_{d-h}^t)\cap I^n=(x_1^t,\ldots,x_{d-h}^t)I^n$ for
 all $n\ge 0$; 
see e.~g.~\cite[Theorem 13.10 and Theorem 13.7]{HIO}).

Furthermore, because $\overline G$ is obtained from $G$ by killing 
elements of degree zero, the $a$-invariants of $G$ and $\overline G$
are equal (see Remark 
\ref{hyperfact} here, or  
~\cite[Remark 5.1.21]{Va}).  Because the Rees ring $A[It]$ is Cohen-Macaulay,
the $a$-invariants $a(G) $ and $a(\overline G)$ are  negative.
So because $\overline G$ is Cohen-Macaulay with
negative $a$-invariant, the corresponding Rees ring
$$
\overline R = \overline A \oplus \overline I \oplus {\overline I}^2 \dots
$$
is Cohen-Macaulay. So Theorem \ref{mainalg} can be applied to 
the $\overline \fm$-primary  ideal $\overline I$ in the Cohen-Macaulay ring 
$\overline A$. 

Now observe that every minimal reduction of $\overline I$ is of 
type can be obtained as the image $\overline J$ in $\overline A$ of
some minimal reduction $J$ of $I$. Indeed,
 let $(\overline{a_1},\ldots,\overline{a_h})$
be a minimal reduction of $\overline I$ with $a_1,\dots,a_h\in I$. 
Then ${\overline I}^{n+1}=
(\overline{a_1},\ldots,\overline{a_h}) I^n$ for some $n\ge 0$ implies 
\begin{align*}
\begin{split}
I^{n+1}&\subset (a_1,\ldots,a_h)I^n+(x_1^t,\ldots,x_{d-h}^t)\cap I^{n+1}\\
&=(a_1,\ldots,a_h)I^n+(x_1^t,\ldots,x_{d-h}^t)I^{n+1}
\end{split}
\end{align*} so that $I^{n+1}\subset (a_1,\ldots,a_h)I^n$.
It thus follows that $(a_1,\ldots,a_h)$ is a minimal reduction of $I$. 

Finally, 
let $y\in \core(I\omega_A)$.
Since $\omega_{\overline A} = 
\frac{\omega_A}{(x_1^t,\ldots,x_{d-h}^t)\omega_A}$, the
 above computation shows that $\overline y \in \core 
(\overline I \omega_{\overline A})$. 
Set 
 $\overline Y=\Proj \overline R$. Then by Theorem \ref{mainalg} together
with Proposition \ref{intersectionlem1}, we have
$$ \core{(\overline I \omega_{\overline A})} \subset  
\Gamma(\overline{Y},\omega_{\overline Y}(h)).$$
Now, again using the injection 
$\Gamma(\overline{Y},\omega_{\overline Y}(n)) \hookrightarrow 
\omega_{\overline A}$
induced by the trace map, one easily checks that 
there is an induced isomorphism
$$
\Gamma(\overline{Y},\omega_{\overline Y}(n)) \cong 
\frac{\Omega_{n}}{(x_1^t,\ldots,x_{d-h}^t)\Omega_n}
$$ for all 
$n\in \bbZ$.
Thus 
$$ \overline y := y {\text{ mod }} (x_1^t,\ldots,x_{d-h}^t)\omega_A 
\in \frac{\Omega_{h}}{(x_1^t,\ldots,x_{d-h}^t)\Omega_h} \subset 
\frac{\omega_{A}}{(x_1^t,\ldots,x_{d-h}^t)\omega_A},
$$
so 
$$ y \in \Omega_h + (x_1^t,\ldots,x_{d-h}^t)\omega_A.
$$
Finally, because this works for any positive $t$, we have 
  \begin{equation*}y\in 
\bigcap_{t\ge 1}(\Omega_h+(x_1^t,\ldots,x_{d-h}^t)\omega_A)
=\Omega_h.\end{equation*}
This shows that $\core(I\omega_A)\subset \Omega_h$, and the proof is complete.
\end{proof}

\subsection{Core in Dimension one.}

The main technical theorem easily gives a formula for the
core of an $\fm$-primary ideal in a local ring of dimension one, without
any Cohen-Macaulay hypothesis at all. 

\begin{cor}
Let $(A, \fm)$ be a one dimensional local ring 
containing  the rational numbers.
Then for any $\fm$-primary ideal $I$, we have 
$$ \core(I\omega_A) =  \Omega_1.$$
In particular, if $A$ is Cohen-Macaulay, then 
$$
\core(I) = J^{r+1}: I^r,
$$
where $J$ is any reduction of $I$ and $r$ is any positive integer
such that $I^{r+1} = JI^r$. 
\end{cor}

\begin{proof}
The second statement follows immediately from the first using 
Lemma \ref{lem6}. The first statement follows immediately from the 
 the main technical theorem, Theorem \ref{mainalg}.
One need only verify that $\Omega_1 \subset J\omega_A$ and
 $J\omega_A \cap \Omega_0 = 
J(\Omega_{-1} \cap \omega_A)$, but this is trivial by 
Lemma \ref{lem1} since $J$ is generated by a non-zero-divisor on 
$\omega_A$.
\end{proof}

\subsection{Core and  adjoints.} \label{adjoint}
We recall the definition of an adjoint  (or multiplier) ideal.
Although a definition can be given that does not refer to
resolution of singularities (see \cite{L3}), we prefer the 
following approach.

Let $X$ be a Gorenstein scheme 
 essentially of finite type over a field
of characteristic zero, and let $\mathfrak {a}$ be a
coherent sheaf of ideals on $X$. Fix a log resolution of $\fa$,
that is, a  proper birational map $Y \overset{\pi}\to X$ from a 
 smooth scheme  $Y$ such that $\fa\cO_Y$ is locally principal 
and the union of the  
support of the 
corresponding divisor and the exceptional divisors is a divisor with 
normal crossing support. 
Then the multiplier (or adjoint) ideal of $\fa$ 
is the ideal sheaf of $\cO_X$
$$
\adj(\fa) = \pi_*(\fa \omega_{Y/X}),
$$
where $\omega_{Y/X} = \omega_Y\otimes \pi^*\omega_X$
 is the relative canonical sheaf of $\pi$.
This is independent of the choice of the log resolution. 
Note that because $X$ is Gorenstein, $\omega_X$ is invertible,
so $\omega_{Y/X}$ is invertible as well. 
See  \cite{laz} or  \cite{Ein} for the general theory of multiplier ideals
from the algebro-geometric point of view, or \cite{L3} for a more
algebraic point of view.

\bigskip
In \cite{HS}, 
 Huneke and Swanson studied 
 the core of an integrally closed $\fm$-primary 
ideal in a two-dimensional regular local ring. 
In particular, they showed that in this case,
$$
\core{(I)} =  \adj{(I^2)} = I \adj{(I)}.
$$
However, such ideals are very special in a sense: 
 the corresponding Rees algebra always has rational singularities 
(see~\cite[Proposition 1.2]{L1} and~\cite[Proposition 2.1]{Hy1}).
 In particular, it is Cohen-Macaulay and normal.
 For this reason the following 
 can be considered a natural generalization of their result to higher
dimension.

\begin{cor} 
\label{ratsingcor}
Let $A$ be regular 
local ring essentially of finite type over a field of characteristic zero. 
Let $I\subset A$ be equimultiple ideal of positive height 
$h$ such that the Rees ring $A[It]$ is normal and Cohen-Macaulay. Then the following conditions
are equivalent
\begin{itemize}
\item[1)] $A[It]$ has rational singularities;
\smallskip
\item[2)] $\Omega_n=\adj(I^n)$ for all $n\ge 0$;
\smallskip
\item[3)] $\core(I)=\adj(I^h)$.
\end{itemize}
If this is the case, then 
$$
\core(I)= I \adj(I^{h-1}) {\text{  and  }}
 \adj(I^{h-1})=\core(I):I.$$ 
\end{cor}

\begin{rem}
In fact,  as is clear from the
proof, it is not necessary that $A$ be regular for this theorem.
It is sufficient if $A$ is 
Gorenstein with rational singularities.
\end{rem}

\begin{rem}
In dimension two it is 
 not necessary that $A$ is essentially of finite type over a field of characteristic zero. Indeed, resolutions exist in this setting, and 
the hypotheses on the Rees ring 
 imply that the reduction number is at most one;
see Remark \ref{304}. 
In particular, since the hypotheses on the Rees ring 
hold automatically for any
$\fm$-primary integrally closed ideal in a two regular local
 dimensional local ring, 
 the  Huneke-Swanson theorem is recovered in full generality.

\end{rem}

\begin{proof} Set $Y=\Proj A[It]$.  Observe first that $A[It]$ has rational singularities if and only if $Y$ has rational singularities
(see for example, ~\cite[Proposition 1.2]{L4} and~\cite[Proposition 2.1]{Hy1}).
Let $f\colon Z\lrarrow Y$ be 
a log resolution of $I\cO_Y$.
Because $Y$ is Cohen-Macaulay and normal, it follows that 
 $Y$ has rational singularities if and only
if the natural inclusion 
$$f_*\omega_Z \subset \omega_Y$$ is an isomorphism;
\cite{Kempf}.
On the other hand, 
because $I\cO_Y$ is ample for the map $Y \la \spec A$, 
 this equivalent to requiring that the natural map
$$\Gamma(Y,I^n\cO_Y \otimes f_*\omega_Z) \hookrightarrow 
\Gamma(Y,I^n\omega_Y)$$
be an isomorphism for all
$n\gg 0$. Because 
$I^n\cO_Y \otimes f_*\omega_Z$ can be identified with 
$ f_*(I^n\omega_Z)$, this is the same as the
the natural inclusion 
\begin{equation}\label{eq7}
\adj(I^n) = 
\Gamma(Z, I^n\omega_Z) \hookrightarrow \Gamma (Y, I^n\omega_Y) = \Omega_n
\end{equation}
being an isomorphism for all $n \gg 0$.

Now, 
by Lipman's Brian{\c c}on-Skoda Theorem (see also \cite{laz})
$$\adj(I^h) = I\adj(I^{h-1})$$
and by our Brian\c con-Skoda Theorem (actually Corollary \ref{BS2} and 
the subsequent remark)
$$
\Omega_h = I\Omega_{h-1}.$$
So remembering also that $\Omega_{n+1}:_{\omega_A}I
 = \Omega_n$ for all $n \geq 0$, by Lemma \ref{colonlem}, we conclude that 
  (\ref{eq7}) is an isomorphism
  for all $n \gg 0$ if and only if it is an isomorphism for 
 $n  = h-1$. This proves the equivalence of statements (1) and (2).
The equivalence with (3) is also clear, since $\core(I) = \Omega_h$
by Theorem \ref{corethm}. 

Finally, 
the formula $\adj(I^{h-1})=\core(I):I$ is a consequence of the formula
$\Omega_{h-1} = \Omega_h:_A I$ of Lemma  \ref{colonlem}. The corollary
is proved.
\end{proof}

\begin{example}
Let $I$ be a normal equimultiple
 monomial ideal of height $h$  in a polynomial ring $S$
over $\CC$.
Then 
$$
\core I   = \adj(I^h) = I\adj(I^{h-1}).
$$
Indeed, 
in this case, the Rees ring $S[It]$
is a normal semi-group algebra, and hence has rational singularities
(since it is a direct summand of a polynomial ring
\cite{bou}).
In particular, if $I$ is generated by monomials 
$x^{A_i} = x_1^{a_{i1}}x_{2}^{a_{i2}}\dots x_n^{a_{id}},
$ 
then  $\core(I)$ is generated by monomials $x^B = x_1^{b_1} \dots x_d^{b_d}$
where
$$
(b_1 + 1, b_2 + 1, \dots, b_d + 1) 
$$
is in the interior of the
 convex hull of the points $h A_1, \dots, h A_r$ in $\NN^{d}$; 
\cite{how}.
\end{example}

\subsection{Core and coefficient ideals.}
If $J\subset I$ is a reduction of $I$, Aberbach and Huneke defined the
{\it  coefficient ideal} $\fa(I,J)$ as the largest 
ideal $\fa$ such that
$I\fa=J\fa$; see \cite{AH}. 
The next corollary relates this notion to the core.

\begin{cor} \label{coeffcorecor}Let $(A,\fm)$ be a Gorenstein local ring containing the set of rational
numbers.  Let $I\subset A$ be an equimultiple ideal of positive height
$h$ such that the Rees ring  $A[It]$ is Cohen-Macaulay. If 
$J\subset I$ is any minimal reduction, then 
\begin{equation*}
\core(I)=I\fa(I,J)\quad\text{and}\quad\fa(I,J)=\core(I):I.\end{equation*}
\end{cor}

\begin{proof} By Theorem~\ref{corethm} we know that $\core(I)=\Omega_h$. By the Brian\c con-Skoda Theorem 
(Lemma~\ref{BS}), $\Omega_h=I\Omega_{h-1}$. On the other hand, according
to~\cite[Theorem 3.4]{Hy2},  
$\Omega_{h-1}=\fa(I,J)$. Thus the first claim follows. 
The second one is now a consequence of the formula 
$\Omega_{h-1}=\Omega_h:I$ of Lemma~\ref{colonlem}. \end{proof}

\begin{rem} 
In fact, rational 
singularities of $A[It]$  can be characterized in terms of the equality
$\fa(I,J)=\adj(I^{h-1})$ where $J\subset I$ 
is any minimal reduction; See \cite[Corollary 3.5]{Hy2}. 
\end{rem}

\subsection{Further Properties of Core, and Questions}

\begin{cor} Let $A$ be a Gorenstein local ring
essentially of finite type over a field of characteristic zero. Let
$I$ be  an equimultiple ideal of positive height $h$ such 
that the Rees ring
$A[It]$
 has  rational singularities. Then 
$$\core(I) \subset \core(I')$$
for any ideal $I'$ of height $h$ containing $I$. 
\end{cor} 

\begin{proof} 
Because  $I\subset I'$, we know $\adj(I^h) \subset \adj{((I')^h})$.
By  Corollary~\ref{ratsingcor}, we have $\core(I) = \adj(I^h) \subset 
\adj{((I')^h)}$. On the other hand, 
from the Brian\c con-Skoda theorem (\cite{L3}), 
$\adj{((I')^h)}$ is contained in every reduction of $I'$.
So $\core(I) \subset \core(I')$.
 \end{proof}

\begin{question}\label{q1}
If $I$ is an 
integrally closed ideal,
then 
is 
 $\core(I) \subset
\core(I')$
for all ideals $I'$ containing $I$? 
If $I$ is not integrally closed, the answer is no in general.
Indeed, whenever
 $I \subset I'$ is an integral extension of ideals,
then a minimal reduction of $I$ is a minimal reduction of  $I'$,
but $I'$ may admit reductions that are not reductions of $I$.
So clearly 
 $\core(I') \subset \core(I)$, but the inclusion can be strict,
and usually is, for example, when 
$I $ is a minimal  reduction of $I'$.
On the other hand, the same reasoning
indicates that 
there is no loss of generality in assuming that  
 also that  $I'$ is  
integrally closed in Question \ref{q1}.  
This question was first raised in \cite{HS}.
\end{question}

The next result has to do with when the core itself is
integrally closed. 

\begin{prop}\label{icprop}Let $A$ be a Cohen-Macaulay 
local domain of dimension $d$ containing the set of 
rational numbers. Let $I\subset A$ be a normal
equimultiple ideal of positive
height $h$ such that  the Rees ring 
$A[It]$ is Cohen-Macaulay. 
 Then $\core(I)$ is an integrally closed ideal of $A$.
\end{prop}

\begin{proof} Set $Y=\Proj A[It]$. Consider $\omega_Y$ as
 a sub-sheaf of the constant sheaf $K$ 
where $K$ is the quotient field of $A$. 
Since $\omega_Y$ is reflexive, we know that
\begin{equation*}\Omega_h=\Gamma(Y,I^h\omega_Y)=
\bigcap_{\codim\overline{\{x\}}=1} I^h\omega_{Y,x}.\end{equation*}
 But
then 
\begin{equation*}
\core(I)=\Omega_h:_A\omega_A=\bigcap_{\codim\overline{\{x\}}=1}
I^h\omega_{Y,x}:_A \omega_A.\end{equation*}
Thus $\core{(I)}$ 
 is integrally closed, because it is an 
intersection of
 integrally closed ideals of
$A$. \end{proof}

\begin{question}
Under what conditions is the 
 core of a normal ideal  integrally closed? 
This issue was first raised in \cite{HS}. See also
  \cite[Examples 3.9 and 3.10]{CPU2}.
\end{question}
\bigskip

Finally, we record an observation about the asymptotic behavior of 
core, as a partial answer to a question raised in \cite{HS}.

\begin{cor} Let $(A,\fm)$ be
 a Gorenstein local ring containing the rational numbers.
Let $I\subset A$ be an equimultiple ideal of positive height $h$
 such that $A[It]$ is Cohen-Macaulay. Then
\begin{equation*}\core(I^n)=I^{(n-1)h}\core(I)\end{equation*}
for all $n\ge 1$. 
\end{cor}

\begin{proof} Set $Y=\Proj A[I^nt]$.
 Since $Y=\Proj R^{(n)}\cong \Proj R$, we observe that
$\Gamma(Y,I^{nk}\omega_Y)=\Omega_{kn}$ for all $k\in \bbZ$. The ideal 
$I^n$ being also
equimultiple of height $h$,
 Theorem~\ref{corethm} and Corollary \ref{BS2} (and the subsequent remark),  
 now give
\begin{equation*}\core(I^n)=\Gamma(Y, I^{nh}\omega_Y)
=\Omega_{nh}=I^{nh-h}\Omega_h=
I^{(n-1)h}\core(I).\end{equation*}
\end{proof}

\section{Non-Vanishing Sections and the Core}\label{geom-to-alg}

The goal of this section is to reduce
 Kawamata's 
 Conjecture
to a purely algebraic statement relating the
core of an $\fm$-primary ideal in a local ring of dimension $d$ to 
the adjoint ideal (or multiplier ideal) of the $d$-th power of the ideal. 
Actually, of course, we must work in the graded category. Also, to 
get at the most general version of Conjecture \ref{mainkaw},
we must expand the notions of the core and the adjoint to
submodules of the canonical module.
The main result of this section is the following theorem.

\begin{thm}\label{maingrcore}
Let $D$ be an ample Cartier divisor on a rationally singular projective 
variety $X$ of positive dimension, and let 
$$
S = \bigoplus_{n \geq 0}  H^0(X, nD)
$$
be the corresponding section ring. 
Fix $n \gg 0$, and let $I = S_{\geq n}$
 be the ideal generated by all elements of degrees at least $n$ in $S$.
Then $H^0(X, D) \neq 0$ if 
\begin{equation}\label{grad}
\grcore{ (I\omega_S)}  \, = \adj (I^{d+1}\omega_S), \,\,\,\,
{\text{ where }}d+1 = \dim S,
\end{equation}
as subsets of $\omega_S$.
\end{thm}

\begin{rem}\label{mult.ideal}
Here, 
$\grcore{( I\omega_{S})}$
 denotes the intersection of all submodules of $\omega_S$
of the form $J\omega_S$, where $J$ is a homogeneous reduction of $I$.
Likewise, for any ideal $I$ in a normal domain $S$, $\adj(I\omega_S)$
denotes the following natural variant of the adjoint ideal. Fix
a log resolution $Y \rightarrow \spec S$ of $I$. Then $\adj(I\omega_S)$
is the submodule of $\omega_S$ given by $\pi_*(I\omega_S)$.
This definition is independent
of the choice of log resolution. This is
like the usual notion of multiplier ideal, but  the
relative canonical modules has been replaced by the absolute
 canonical module of $Y$. This has the advantage of being
defined even when $S$ is not Gorenstein (or $\QQ$-Gorenstein).
However, it is a submodule of $\omega_S$ rather than an ideal of $S$.
\end{rem}

\begin{rem}
Formula (\ref{grad}) is a graded version of the formula we
proved (under certain conditions on the ring and ideal) in 
Section 5. As we will see, rings arising from divisors satisfying the
hypothesis of Conjecture \ref{mainkaw} satisfy these conditions,
so Kawamata's Conjecture is very closely related to our formulas
in \ref{ratsingcor}.
\end{rem}

\begin{rem} As will be clear from the proof, 
a version of Theorem \ref{maingrcore} holds if $X$ is not
necessarily rationally singular, but is only normal. In this case,
$H^0(X, D) $ is non-zero  if $\grcore(I\omega_S) = \Omega_{d+1}$,
with notation as in Section 3. 
\end{rem}

\begin{rem}
In fact, the converse of Theorem \ref{maingrcore} is also true:
$H^0(X, D)$ is non-zero if and only if the formula (\ref{grad}) holds
for $n\gg 0$ in the section ring $S$ of $D$. However, the
proof of the requires rather different ideas and techniques, so
we postpone it to a subsequent paper. 
\end{rem}

In this section, we 
first  
prove 
 Theorem \ref{maingrcore}. We then investigate the 
hypothesis forced upon the section ring $S$ of a pair
$(X, D)$ satisfying the  hypothesis of Conjecture \ref{mainkaw}.
 Finally, we  end with a discussion of
core versus graded core in a graded ring.

\subsection{A general criterion for non-vanishing.} \label{nonvan}

Let $D$ be an ample Cartier divisor on a normal projective variety $X$
 of dimension $ d\geq 1$.
By definition, the section ring of the pair $(X, D)$ is the
$\bbN$-graded ring
$$
S = \bigoplus_{n \in \bbN} H^0(X, \cO_X(nD)),
$$
whose multiplication is given by the natural multiplication of sections.
The ring $S$ is a normal graded domain, finitely generated over the
field $k = H^0(X, \cO_X)$, which we will assume to be infinite.
 There is a canonical isomorphism from $X$ to 
$\proj S$ under which the invertible sheaf $\cO_X(nD)$ corresponds
to the coherent module on $\proj S$ arising from the graded $S$-module 
$S(n)$, where $S(n)$ denotes the module $S$ with its grading shifted so that 
$S(n)_{m} = S_{m+n}$. For a general reference on section rings,
see \cite[Section 4.5]{EGAII}.

\begin{prop}\label{van-alg2}
Let $D$ be an ample divisor on a normal 
 projective variety $X$, and let $S$ be the corresponding
section ring of the pair $(X, D)$. 
Fix $n \gg 0$.
Then 
$$
H^0(X, D) = 0$$
if and only if 
$$
[\omega_S]_{n(d+1) - 1} \subset \grcore (I\omega_S),
$$
where $\omega_S$ is the (graded) canonical module of the normal
ring $S$, and $I$ is the ideal of $S$ generated by elements of degrees
at least $n$. 
\end{prop}

Proposition \ref{van-alg2} follows readily from the
following very general criterion for the 
vanishing of the space of global sections of an ample line bundle.

    \begin{lem}\label{van-alg}
Let $D$ be an ample Cartier divisor on a normal 
 projective variety $X$ of
dimension $ d \geq 1$. Fix any integer $i$. Then
$H^0(X, iD) = 0$ if and only if for some (equivalently, every) $n \gg 0$,
and any set 
 $x_0, \dots, x_d$ of $d+1$ global generators for $\cO_X(nD)$,  
the natural inclusion
$$
\sum_{i=0}^{d} x_iH^0(X, K_X + (nd - i)D)
\subset 
H^0(X, K_X + [n(d + 1)  - i]D  )
$$
is an equality.
{\footnote
{The precise condition on  $n$ in Propositions 
\ref{van-alg2} and \ref{van-alg}
 is that $n$ should be large enough that 
$\cO_X(nD)$ is
 globally generated and $H^i(X, K_X + mD) $ vanishes
 for
all $i > 0$ and all $m \geq n-i$.}}
\end{lem}

\begin{proof}[Proof of Lemma \ref{van-alg}]
Fix any set of $d+1$ global generators for $\cO_X(nD)$.
Such a set always exists (assuming $X$ to be defined over an infinite field),
because we can  take generic linear combinations of any set of 
global generators for  $\cO_X(nD)$.

Consider the Koszul complex determined by the $x_i$'s:
$$
0 \la \cO_X(-(d+1)nD) \la \dots \la \bigoplus_{i=0}^d \cO_X(-nD) 
\overset{\begin{pmatrix}x_0 \\ \vdots \\ x_d\end{pmatrix}}{\la} \cO_X \la 0.
$$
Because the $x_i$'s generate $\cO_X(nD)$, this complex is exact.
Tensoring with the invertible sheaf $\cO_X(K_X + [(d+1)n -i]D)$,
we get an exact complex
$$
0 \la
 \cO_X(K_X - i D)
 \la 
\dots
 \la
 \bigoplus_{i=0}^d \cO_X(K_X + (dn - i)D) 
\overset{\begin{pmatrix}x_0 \\ \vdots \\ x_d\end{pmatrix}}{\la} 
\cO_X(K_X +[(d+1)n -i ]D) \la 0.
$$
Because $H^i(X,  \cO_X(K_X + mD)) = 0$ for all $m \geq n-i$ and all $i >0$,
a standard  argument
{\footnote{The standard argument is this: break the complex into several short
exact complexes of sheaves. Then look at the corresponding long 
exact complexes of cohomology, beginning with the $0$-th cohomology of
the short exact sequence arising from the right-most part of the complex.
Working backwards,  the relevant cohomology is the $i$-th
cohomology of the $i$-th short exact sequence from the right.
A similar argument is written down in full in the proof of 
Proposition \ref{intersectionlem1}.
}}   
 shows that
the map of global sections
$$
\bigoplus_{i+0}^d H^0(X,  K_X + (dn - i)D)  
\overset{\begin{pmatrix}
x_0 \\ \vdots \\ x_d\end{pmatrix}}{\la} H^0(X, K_X +[(d+1)n -i]D) 
$$
is surjective if and only if $ H^d(X, K_X -i D) = 0$.

By Serre duality (which holds at the ``top spot'' even if $X$
 is not Cohen-Macaulay),
the claim that $H^0(X, i D)$ is  zero is identical to the claim
that $H^d(X, K_X - i D)$ is zero.
So
$H^0(X, D)$ vanishes  if and only if
$$
H^0(X,  K_X + [n(d+1) -i ]D) \subset (x_0, \dots, x_d)
H^0(X, K_X + [nd - i]D).
$$ 
The proof of Lemma \ref{van-alg} is complete.
\end{proof}

\begin{proof}[Proof of Proposition \ref{van-alg2}]
Interpret global sections of $\cO_X(nD)$ as degree $n$ elements of
$S$. Because   a set of global sections $\{x_i\}$ generates
$\cO_X(nD)$ if and only if their common zero set on $X = \proj S$ is empty,
such a set  is a generating set for 
$\cO_X(nD)$ if and only if the elements $\{x_i\}$ in $S$ generate
 an $\fm$-primary
ideal of $S$. In particular, a set of $d+1$ global sections of 
$\cO_X(nD)$  is a generating set if and only if its elements
 form a homogeneous 
system of
parameters for the ring $S$.

Now fix $n \gg 0$.  If $H^0(X, D)$ is zero, then for each set 
$\{x_0, x_1, \dots, x_d\}$ of global generators of $\cO_X(nD)$, 
Lemma \ref{van-alg} ensures that 
$$
\sum_{i=0}^d x_iH^0(X, K_X + (nd -1)D) = H^0(X, K_X + (n(d+1) -1)D).
$$ 
Interpreted in terms of the section ring $S$, this says
$$
(x_0, \dots, x_d)[\omega_S]_{nd-1} = [\omega_S]_{n(d+1)-1}.
$$
In particular, 
$$
[\omega_S]_{n(d+1)-1} \subset (x_0, \dots, x_d)\omega_S
$$  
for  {\it every } system of parameters for $S$ consisting of
elements of degree $n$. In other words,
$$
[\omega_S]_{n(d+1)-1} \subset 
\bigcap_{J{\text{ s.o.p degree $n$.}}} J \omega_S,
$$
where $J$ ranges over all homogeneous systems of 
parameters for $S$ consisting of elements of degree $n$. 
By Proposition \ref{homred} 
 a system of parameters of degree $n$ is
precisely the same as a minimal 
homogeneous reduction for the ideal $I = S_{\geq n}$,
the ideal generated by all elements of degrees at least $n$,
so this means that 
$$
[\omega_S]_{n(d+1)-1} \subset \grcore (I\omega_S).
$$
The proof of the converse simply reverses this argument.
\end{proof}

\begin{rem}
Although Proposition \ref{van-alg2} follows quite trivially from 
Lemma \ref{van-alg}, the passage to this more algebraic
lemma seems powerful. The point is that there is hope 
for showing that the 
{\it intersection} over all submodules of the form $(x_0, \dots, x_d)\omega_S$
is quite small, so small in fact, that it 
can not contain any element of degree $n(d+1) -1$  (although each 
individual module $(x_0, \dots, x_d)\omega_S$ certainly contains many such 
elements!). 
This would settle Kawamata's Conjecture if it could be accomplished.
\end{rem}

\subsection{An Adjoint Computation.}\label{adjointcomp}

The next proposition is a very general  computation of
 adjoint modules for certain types of ideals in a section ring.
\begin{prop}\label{compute}
Let $S$ be a section ring of a pair $(X, D)$ consisting of an
ample Cartier divisor on a normal projective variety.
 Fix $n \gg 0$ and let $I = S_{\geq n}$
 be the ideal of $S$ generated by all elements of degrees
at least $n$. 
Then 
$$
\Gamma (Y, I\omega_Y) = [\omega_S]_{\geq n+1}, 
$$
as submodules of $\omega_S$, 
where $Y$ is the blowup of the scheme $\spec S$ along the ideal $I$. 
(The precise condition on $n$ is that $n$ should be large enough that 
$(S_{\geq n})^t = S_{\geq tn}$ for all $t \geq 0$.)
\end{prop}

\medskip
The proof  of Proposition \ref{compute}
makes use of the ``natural'' construction from \cite[Section 8.7.3]{EGAII}.

\subsubsection{The ``natural'' construction.} \label{natural}
Let $S^{\natural}$ be the graded ring
$$
S^{\natural} = S \oplus S_{\geq 1} \oplus S_{\geq 2} \oplus S_{\geq 3} \dots,
$$
where $S_{\geq n}$ indicates the ideal of  $S$ generated by 
elements of degree at least $n$. 
The ring 
 $S^{\natural}$ is finitely generated over its degree zero part $S$,
so for large $n$, we have $(S_{\geq n})^k = S_{\geq nk}$ for all $k \geq 0$.
Since the projective scheme of a graded ring is unchanged
under passing to any Veronese sub-ring, we have
$$
 \proj  S^{\natural} \cong \proj S[It]$$
where $S[It]$ is the Rees ring of $S$ with respect to the ideal
$I = S_{\geq n}$. 

There are two natural geometric interpretations of the scheme 
$\proj S^{\natural}$.
First, the above isomorphism shows that  $\proj S^{\natural}$ 
can be considered 
as the blowup of the ideal $I$ in the affine scheme $\spec S$. 
On the other hand, there is
a natural isomorphism  (see \cite[8.7.3]{EGAII})
$$
\spec_X(\cO_X \oplus \cO_X(D) \oplus \cO_X(2D) \oplus \dots )
\overset{\cong}{\la}
 \proj (S \oplus S_{\geq 1} \oplus S_{\geq 2} \oplus \dots). 
$$
This allows to interpret $  \proj  S^{\natural}
$ also as the total space of the
``tautological'' line bundle $\cO_X(-D)$ on $X$ (or, in some
writers' terminology, as the scheme $\bbV(\cO_X(D))$).
Correspondingly, there are two natural  projections,
$$
 \proj S^{\natural} \overset{\pi}{\la} \spec S 
\qquad \qquad
 \proj S^{\natural} \overset{\eta}{\la} \proj S = X.
$$
The first is the blowing up morphism, while the second is the
structure map of the line bundle $\cO_X(-D)$. 

\smallskip
\begin{rem} 
These interpretations of $\proj 
S^{\natural}$   generalize the following
situation. Let $Y$ denote the incidence correspondence
$$
Y = \{(p, \ell) \, | \, p \in \ell\} \subset \bbC^n \times \bbP^{n-1},
$$
where $\ell$ is a line through the origin in $ \bbC^{n}$
and $p$ is a point on it. By projecting $Y$ to 
either $\bbC^n$ or $\bbP^{n-1}$, respectively, we arrive
at either the blowup of the origin in $\bbC^n$ or the 
structure map of the tautological line bundle on $\bbC^n$.
Note that here $Y = \proj S^{\natural}$ where 
$S$ is the polynomial ring in $n$ variables. 
\end{rem}

\begin{proof}[Proof of Proposition \ref{compute}]
Using the interpretation of $Y$ as the total space of the tautological
bundle, one easily computes  $\omega_Y$. Indeed, because
$Y \overset{\eta}{\to} X$ is smooth of relative dimension one with 
$\Omega_{Y/X} = \eta^*\cO_X(D)$, we have
$$
\omega_Y = \eta^*\omega_X \otimes \Omega_{Y/X} = \eta^*\omega_X \otimes
\eta^*\cO_X(D) = \eta^*(\omega_X(D)).
$$ 
Also, thinking of $S[It]$ as the $n$-th Veronese sub-ring of the
algebra $S^{\natural}$, we see that $I^t\cO_Y = \cO_Y(nt) = \eta^*\cO_X(ntD)$,
where $\cO_Y(nt)$ is the coherent sheaf on $Y = \proj S^{\natural}$
corresponding to the graded module $S^{\natural}(nt)$.

Now, noting that the map $Y \overset{\eta}{\to} X$ affine
and that $\eta_*\cO_Y = \bigoplus_{i \in \bbN} \cO_X(iD)$,  we compute
\begin{align*}
\Gamma(Y, I^{t}\omega_Y) = &
\Gamma (Y, \eta^*\cO_X(tnD) \otimes \eta^*(\omega_X \otimes \cO_X(D)) ) \\
= & \Gamma (Y, \eta^*(\omega_X \otimes \cO_X([tn+1]D)) ) \\
= &  \Gamma (X, (\eta_*\cO_Y \otimes (\omega_X \otimes \cO_X([tn+1]D)) ) \\
= &  \Gamma (X, (\bigoplus_{i\in \bbN}\cO_X(iD)) 
 \otimes (\omega_X \otimes \cO_X([tn+1]D)) ) \\
=&  \bigoplus_{i\in \bbN} \Gamma(X, \omega_X([tn+1+i]D)\\
= & [\omega_S]_{\geq tn+1}.
\end{align*}
This completes the proof of Proposition \ref{compute}.
\end{proof}

\smallskip
\begin{proof}[Proof of Theorem \ref{maingrcore}]
Suppose that $H^0(X, D) = 0$. Then by Proposition \ref{van-alg2},
we have that 
$$ 
[\omega_S]_{n(d+1) -1} \subset \grcore(I\omega_S).
$$
If the equality (\ref{grad}) of Theorem \ref{maingrcore}  holds,
then in fact
$$   
[\omega_S]_{n(d+1) -1} \subset \adj(I^{d+1}\omega_S).
$$

Note that the assumptions of Kawamata's Conjecture force the variety $X$
to have rational singularities. Therefore, 
 the scheme $Y = \proj S[It]$ also
 has rational singularities, because it is the
total space of a line bundle over $X$. Furthermore, the 
ideal $I$ (and its powers) pull back to an invertible sheaf on $Y$ under
the birational map $Y \la \spec S$ whose support is an irreducible
closed subvariety (isomorphic to $X$). It is easy to check that
in this situation, the adjoint of $I^{d+1}$ can be computed from the
resolution $Y \la \spec S$. If particular, if $H^0(X, D) = 0$, then
 $$   
[\omega_S]_{n(d+1) -1} \subset \Gamma(I^{d+1}\omega_Y).
$$
Finally, since $I^{d+1} = S_{\geq n(d+1)}$,
we apply Proposition \ref{compute} to conclude that
$$   
[\omega_S]_{n(d+1) -1} \subset [\omega_S]_{\geq n(d+1)+1}.
$$
This is an obvious contradiction, since $\omega_S$ is non-zero in
all sufficiently large degrees. Thus $H^0(X, D)$ can not vanish
and the proof is complete.
\end{proof}
\smallskip
\subsection{The section ring.}\label{sectionring}
In order to use Theorem \ref{maingrcore} to 
prove Kawamata's Conjecture,  we need to better understand
the special conditions imposed on $S$  by the 
hypothesis of Conjecture \ref{mainkaw}. Remarkably, 
it turns out that just the right condition to deduce a local form 
of formula (\ref{grad}) is satisfied.

\begin{prop}\label{CMirr}
Let $S$ be the section ring of a normal projective variety
 $X$ with respect to
an ample Cartier divisor $D$, and let $I$ be the ideal of
$S$ generated by elements of degree $n\gg 0$. 
Assume that there exists an effective $\bbQ$-divisor $B$ 
such that the pair $(X, B)$ is Kawamata log terminal and 
the $\bbQ$-divisor $D - (K_X + B)$ is big and nef.
Then the irrelevant ideal 
$$
(S[It])^+ : = I \oplus I^2 \oplus I^3 \oplus \dots 
$$
of the Rees ring $S[It]$
is a Cohen-Macaulay $S[It]$-module.
\end{prop}

The proof makes use of the following two lemmas.

\begin{lem}\label{nice}
Let $S$ be a section ring  as in Proposition \ref{CMirr}.
Then the local cohomology modules of $S$ with support in the unique homogeneous maximal ideal $\fm$ satisfy:
\begin{enumerate}
\item For $i < \dim S$, the graded $S$-modules $H^i_{\fm}(S)$
vanish in all degrees $n \neq 0$.  
\item $H^{\dim S}_{\fm}(S)$ vanishes in positive degrees.
\end{enumerate}
\end{lem}

\begin{lem}\label{R+CMchar}
Let $(A, \fm)$ be an arbitrary local 
 ring of dimension $d$
 and let $\{I_n\}$ be a
Noetherian filtration of $A$ consisting of ideals of positive height.
Let 
$$
R : = A \oplus I_1 \oplus I_2 \oplus \dots 
$$
and 
$$
G : = A/I_1 \oplus I_1/I_2 \oplus I_2/I_3 \oplus \dots
$$
denote, respectively, the Rees ring and the associated graded ring of 
$A$ with respect to this filtration, whose unique homogeneous
maximal ideals will be denoted $\fM_R$ and $\fM_G$, respectively. 
Then the irrelevant ideal of $R$ is 
a Cohen-Macaulay $R$-module if and only if 
the following two conditions on $G$ are satisfied:
\begin{enumerate}
\item For $i < d$,  $[H^i_{\fM_G}(G)]_{n} = 0 $ for $n \neq 0$.
\item $H^d_{\fM_G}(G)  $  vanishes in positive degrees.
\end{enumerate}
\end{lem}

\begin{proof}[Proof of Proposition \ref{CMirr}]
Fix $n \gg 0$. 
Setting $I = S_{\geq n}$ to be the ideal  in $S$ generated by all 
elements of degrees at least $n$, we have 
 $I^k = S_{\geq nk}$. 
So for this $n$, 
the Rees ring $S[It]$ is the $n$-th Veronese subring 
of the Rees ring $S^{\natural}$ formed from the filtration $I_n = S_{\geq n}.$
The irrelevant ideal of $S[It]$ is therefore the $n$-th Veronese
submodule of the irrelevant ideal of $S^{\natural}$. Thus in order
to show that the irrelevant ideal of $S[It]$ is a Cohen-Macaulay
 $S[It]$-module, it is sufficient to prove that the irrelevant ideal 
of $S^{\natural}$ is Cohen-Macaulay (since the
appropriate local cohomology modules for the irrelevant ideal of 
$S[It]$ are 
Veronese submodules of the corresponding local cohomology modules
for the irrelevant ideal of $S^{\natural}$).

To show that the irrelevant ideal of $S^{\natural} $
is  Cohen-Macaulay, note that 
 because the irrelevant ideal is graded, it is enough to
check Cohen-Macaulayness after localizing at the unique homogeneous
maximal ideal of $S[It]$. So we may replace $S$ by its localization $A$ at its
unique homogeneous ideal, and replace the filtration by  its image $\{I_n\}$ in
 $A$.
 Note that the 
associated graded ring of $A$ with respect to  this filtration $\{I_n\}$
 is canonically isomorphic
to $S$.
Thus, the Proposition follows immediately from combining 
the two lemmas.
\end{proof}

\begin{proof}[Proof of Lemma \ref{nice}]
First note that because $(X, B)$ is Kawamata log terminal,
  the variety $X$  
 has rational singularities (see e.g. \cite[Th 1.3.6]{KMM}).
 In particular, $X$ is Cohen-Macaulay
and Serre duality  holds for $X$. 

To
check statements (1) and (2), we make use of the identifications
$$
[H^i_{\fm}(S)]_n \cong H^{i-1}(X, \cO_X(nD))
$$
for each $i \geq 2$ and all $n \in \bbZ$ (see \ref{loc-glob}).
Now, for $2 \leq i < \dim S$, the vanishing of $H^{i}_{\fm}(S)$ in negative 
degree  follows from the (dual form of the) Kodaira vanishing theorem
applied to the ample divisor $D$. 
For all $i \geq 2$, the vanishing of $H^{i}_{\fm}(S)$ in positive
degree follows from the Kawamata-Viehweg vanishing theorem. 

Since $S$ is normal, the local cohomology modules $H^1_{\fm}(S)$ and
 $H^0_{\fm}(S)$ are zero in any case, so the Proposition is proved.
\end{proof}

\begin{proof}[Proof of Lemma \ref{R+CMchar}]
Let $R^+$ denote the irrelevant ideal of $R$.
We make use of the following two exact sequences
$$
0 \la R^+ \la R \la A \la 0
$$
and 
$$
0 \la R^+(1) \la R \la G \la 0.
$$
Because the $R$-module $A$ is concentrated in degree zero, 
the first sequence gives 
$$
[H^i_{\fM_R}(R^+)]_n \cong  [H^i_{\fM_R}(R)]_n 
$$
for all $n \neq 0$ and all $i$. 
Then looking at the long exact sequence arising from the second sequence,
we  find a long exact sequence
$$
\dots \la [H^{i-1}_{\fM_R}(G)]_n \la 
 [H^i_{\fM_R}(R^+)]_{n+1} \la 
[H^i_{\fM_R}(R)]_n \la 
 [H^{i}_{\fM_R}(G)]_n \la \dots 
$$
in each degree $n$.

Assume that $R^+$ is Cohen-Macaulay.
Then for $i < d+1 = \dim R^+$, we have 
$[H^i_{\fM_R}(R^+)]_n = 0 $ for all $n$.
This implies that $  [H^i_{\fM_R}(R)]_n = 0 $ for all $n \neq 0$.
So the long exact sequence above tells us that $H^{i}_{\fM_R}(G)$
vanishes in every non-zero degree, for all $i < d$.

For $ i = d$, if $n > 0$, the long exact sequence above becomes
$$  
\la [H^d_{\fM_R}(R)]_n \la [H^d_{\fM_G}(G)]_n \la 
[H^{d+1}_{\fM_R}(R^+)]_{n+1} \la 0.
$$
Because  $H^{d+1}_{\fM_R}(R^+)]_{n+1} \cong 
[H^{d+1}_{\fM_R}(R)]_{n+1}$, and 
  the $a$-invariant of the Rees ring $R$ is $-1$,
we see $ [H^{d+1}_{\fM_R}(R^+)]_{n+1} = 0$ for $n \geq  0$.
So 
both 
modules 
$$
[H^d_{\fM_R}(R)]_n \, \, {\text{  and  }} \, \, [H^{d+1}_{\fM_R}(R^+)]_{n+1}
$$
are zero for $n > 0$, and so 
$H^d_{\fM_G}(G)$ vanishes in positive degree.
The converse argument just reverses this. The lemma is proved.
\end{proof}

\bigskip
\subsection{Core and Graded Core in Graded Rings}
We have seen  that Kawamata's Conjecture follows from the 
following conjecture.

\begin{conj}
Let $S$ be the section ring of a pair $(X, D)$ satisfying the 
hypothesis of Conjecture \ref{mainkaw}. Then 
\begin{equation}\label{eq**}
\grcore(I\omega_S) = \adj(I^{d} \omega_S)
\end{equation}
where $d$ is the dimension of $S$,
and $I = S_{\geq n}$ is the $\fm$-primary ideal generated by elements of degrees
at least $n$, for some $n\gg 0$. 
\end{conj}

On the other hand, 
for a local ring $(A, \fm)$ and
an ideal $I$ satisfying conditions satisfied by $S_{\geq n}$ in
such a section ring, we have proved (see the remark following
 Theorem \ref{corethm} and
 Proposition \ref{CMirr})
 that 
$$
\core (I\omega_A) = \adj (I^d\omega_A).
$$
In particular, $
\core (I\omega_S) = \adj (I^d\omega_S),
$
for the section ring $S$ of a pair $(X, D)$ satisfying the
hypothesis of Kawamata's Conjecture. 
It is easy to believe that perhaps the core and 
graded core of $I = S_{\geq n}$ are equal for large $n$, and hence that 
we have proved Kawamata's Conjecture. However, the 
problem is appears to be quite subtle. In fact, we have 
the following corollary of the Main Technical Theorem.

\begin{cor}\label{cor1}
Let $S$ be a section ring of a normal Cohen-Macaulay
 variety of characteristic zero
with respect to
any ample divisor.
 Let $I = S_{\geq n}$
 be the ideal generated by the homogeneous elements
of $S$ of degrees at least $n$, for $n \gg 0$.
Then 
\begin{equation}\label{eq8}
\core(I\omega_S) = [\omega_S]_{\geq nd +1},
\end{equation}
where $d = \dim S$.
In particular, if the variety is rationally singular, then
 $$
\core(I\omega_S) = \adj(I^{d}\omega_S).
$$
Furthermore, if $S$ is Cohen-Macaulay, 
$$\core(I) = S_{\geq nd + a +1}$$
where $a$ is the $a$-invariant of $S$. 
\end{cor}

\medskip
\begin{rem}\label{gradnotcore}
Of course, since there are  ample line bundles on smooth varieties
with no sections,
 Proposition \ref{van-alg2} makes 
 clear that  formula (\ref{eq**}) can not hold 
 in general. 
Indeed, 
 given any ample line bundle with 
no non-zero 
global sections, one can generate examples of ideals (namely $S_{\geq n}$
for $n \gg 0$)
 in a 
graded ring (the corresponding section ring) which have 
many homogeneous reductions, but  for which the core is not equal to the
graded core. 
\end{rem}

\begin{proof}
The statement may be checked locally at the unique 
homogeneous maximal ideal of $S$, so we can replace $I$ by its expansion to
$A = S_{\fm}$, and our previous results in the local case apply to $S$.
 As always, we let $Y$ denote the blowup of $\spec S$ along 
$I$, and set $\Omega_t = H^0(Y, I^t\omega_Y)$. 
From Proposition \ref{compute}, we
have
$$
\Omega_t = [\omega_S]_{\geq nt+1}  
$$
for all $t$.  Since $S$ is Cohen-Macaulay on the punctured spectrum,
by taking $n \gg 0$, one sees that the irrelevant ideal of 
$S[It]$ is Cohen-Macaulay. Thus 
$$
J\omega_S \cap \Omega_{d-1} = J\Omega_{d-2}
$$
for all reductions $J$ of $I$, by Proposition \ref{intersectionlem2}. 
As pointed out in Remark \ref{one} following the Main Technical Theorem,
this ensures that 
$$
\core I\omega_S \subset \Omega_d = [\omega_S]_{\geq nd+1}.
$$
This holds without requiring that $S$ be Cohen-Macaulay,  as does the
statement 
$\core(I) \subset S_{\geq nd + a +1}.$
Indeed, if $z \in \core{(I)}$ has degree less than $nd + a + 1$, then
by taking any non-zero element $w$ of ${[\omega_S]}_{-a}$, we would have 
an element $yw \in \core{(I\omega_S)}$ of degree less than $nd+1$, 
a contradiction.

For the reverse inclusion, we need also that $ \Omega_d \subset
J\omega_S$. This follows from the Brian\c con-Skoda Theorem \ref{BS}
because  $S$ 
is Cohen-Macaulay on the punctured spectrum and $n$ is sufficiently
large. Indeed, we need only that $H^i(Y, I^{d-1-i}\omega_Y)$ is zero
for all $i \geq 1$, where as usual $Y = \proj S[It]$.
 For  $i < d-1$, this is essentially Serre vanishing (as $n$ is large).
For $i = d-1$, the required 
vanishing holds by Lemma  \ref{van} below.
For  $i \geq d$, all the cohomology vanishes since 
$Y$ has a cover by $d$ open affine sets.

The corresponding  statement for  ideals
follows as in the proof that Theorem \ref{corethm} implies \ref{hswan}.
\end{proof}

\begin{lem}\label{van}
If $I$ is a normal ideal in a normal local ring $A$ of dimension  $d$ at
least two, then
$H^{d-1}(Y,
\omega_Y) = 0$, where $Y = \proj A[It]$.
\end{lem}

\begin{proof}
Let $Z \overset\nu\to Y$ be a resolution of singularities of $Y$.
The composition $Z \to Y \to \spec A$ is a resolution of singularities
of $\spec A$. 
We have a short exact sequence
$$
0 \to \nu_* \omega_Z \to \omega_Y \to Q \to 0
$$
where $Q$ is supported on some set of codimension at least two.
This gives rise to a long exact sequence, which---because $\dim Q$ is at
 most $d-2$--- gives a 
surjection
$$
H^{d-1}(Y, \nu_*\omega_Z) \to H^{d-1}(Y, \omega_Y).
$$
On the other hand, by the Grauert-Riemenschneider vanishing theorem
$R^p\nu_*\omega_Z = 0$ for $p > 0,$ so the appropriate spectral
sequence degenerates to give an isomorphism
$H^i(Y, \nu_*\omega_Z) \cong H^i(Z, \omega_Z)$.
But again by the 
Grauert-Riemenschneider vanishing theorem, this time applied to the
resolution $Z \to \spec A$, we have  $ H^i(Z, \omega_Z)$ vanishes as well
 for all 
$i> 0$.
So $ H^{d-1}(Y, \nu_*\omega_Z)$ must be zero, and therefore, 
so is its surjective image 
$ H^{d-1}(Y, \omega_Y)$.
\end{proof}

Since every normal standard graded domain is a section ring, we have the
following corollary. 

\begin{cor}\label{grade}
Let $S$ be normal Cohen-Macaulay $\NN$-graded domain  finitely generated
by its degree one elements over a field  of characteristic zero. 
Let $\fm$ denote its unique homogeneous maximal ideal. Then  
for all $n >0 $,
\begin{equation}\label{8}
\core (\fm^n) =  \grcore(\fm^n) =  \fm^{nd + a + 1},
\end{equation}
where $d = \dim S$ and $a$ is the $a$-invariant of $S$.
\end{cor}

\begin{proof}
The formula for core  follows from  the above corollary
 since $S_{\geq n} =
\fm^n$. It remains only verify that the core is the graded core in this
situation. Set $I = \fm^n$.
Looking at the proof of the Key Lemma, we see that 
$\grcore(I\omega_S) \subset \Omega_d = [\omega_{S}]_{\geq dn+1}$
if the corresponding intersection 
$\cap_{x_1^*, \dots, x_d^*}(x_1^*, \dots, x_d^*)\omega_G$ is zero in degree
 $d$. But because $I$ is generated by elements all of the same degree,
this follows from Theorem \ref{hard}.
\end{proof}

\begin{rem}
In Corollary \ref{grade}, one can weaken the assumption that $S$ is
Cohen-Macaulay and require only that $S$ has an isolated non-Cohen-Macaulay 
point. Then formula (\ref{8}) holds for $n \gg 0$.
\end{rem}

In order to prove Kawamata's Conjecture then, 
we must understand when core and graded core are
equal. For an $\fm$-primary ideal  generated by elements
of the same degree in a Cohen-Macaulay graded
ring over an infinite field,
 that  core and graded core are equal follows from 
\cite[Theorem 4.5]{CPU}.
But in general, this appears to be  a subtle question.
From the point of view of solving Kawamata's Conjecture
(and  understanding non-emptiness of linear systems more generally),
this question is of great interest for 
 ideals of the form $I = S_{\geq n}$ in a section ring.
We return to this in a subsequent paper.

\newcommand{\noopsort}[1]{} \newcommand{\printfirst}[2]{#1}
  \newcommand{\singleletter}[1]{#1} \newcommand{\switchargs}[2]{#2#1}
\providecommand{\bysame}{\leavevmode\hbox to3em{\hrulefill}\thinspace}
\providecommand{\MR}{\relax\ifhmode\unskip\space\fi MR }
\providecommand{\MRhref}[2]{%
  \href{http://www.ams.org/mathscinet-getitem?mr=#1}{#2}
}
\providecommand{\href}[2]{#2}


\end{document}